\documentclass[11pt]{amsart}
\usepackage{amsmath,amsfonts}
\usepackage{amsmath,amsfonts,latexsym}
\usepackage{amscd,amssymb,amsmath}
\usepackage{amsfonts}
\usepackage{amsbsy}
\usepackage{epsfig,afterpage}
\usepackage{psfrag}
\usepackage{graphicx}
 \newtheorem{lema}{Lemma}[section]
\newtheorem{teo}[lema]{Theorem}
\newtheorem{propo}[lema]{Proposition}
\newtheorem{rema}[lema]{Remark}
\newtheorem{coro}[lema]{Corollary}
\newtheorem{defi}[lema]{Definition}
\newtheorem{exem}[lema]{Example}
\newtheorem{prop}[lema]{Proposition}
\newtheorem{abc}[lema]{}
\newtheorem{assum}[lema]{Assumption}
\newcommand{\bea}{\begin{eqnarray*}}
\newcommand{\eea}{\end{eqnarray*}}
\newcommand{\zz}[1]{}

\newcommand{\hbx}{\hfill$\Box$}

\newcommand{\bz}{\mathbb{Z}}

\newcommand{\br}{\mathbb R}

\newcommand{\bc}{\mathbb C}

\newcommand{\bt}{\mathbb T}

\newcommand{\bn}{\mathbb N}

\begin{document}
\thanks{Supported  by the  Polish Research Grants NCN 2011/03/B/ST1/04533, and  2015/19/B/ST1/01458}
\title[ Geometrically distinct solutions of  variational problems with $O(N)$-symmetry]{Geometrically distinct solutions   given by   symmetries of  variational problems with  the $O(N)$-symmetry}

 \author{Wac{\l}aw Marzantowicz}
 \address{Faculty of Mathematics and Computer Sci.,\, Adam Mickiewicz
University of Pozna{\'n},\, ul. Umultowska 87,\, 61-614 Pozna{\'n},
Poland} \email{marzan@amu.edu.pl}

 \maketitle

 \begin{abstract}
 For  variational problems with   $O(N)$-symmetry    the existence of several
geometrically distinct
 solutions has been  shown   by use of group theoretic  approach in previous articles. It
 was done by a  crafty  choice of a family $H_i \subset O(N)$
 subgroups such that the fixed point subspaces $E^{H_i} \subset E$ of the action in   a corresponding functional space are linearly
 independent, next restricting the problem to each  $E^{H_i}$ and using
 the Palais symmetry principle. In this work we give a thorough
 explanation of this approach showing a correspondence between the
 equivalence classes of such subgroups, partial orthogonal flags in
 $\mathbb{R}^N$, and unordered partitions of  the number $N$. By showing that spaces of functions invariant with respect to different classes of groups are linearly independent we
 prove
 that the amount of series  of geometrically  distinct solutions obtained in this
 way grows exponentially in $N$, in contrast to logarithmic, and
 linear growths of earlier papers.

    \end{abstract}

     \section{Introduction}

     The purpose of this paper is to describe a group theoretical scheme which  arose in works on $O(N)$-invariant  variational problems as a method to show the existence
of several geometrically different series of solutions distinguished
by their  symmetry properties. This approach  originated from the
studies of  problem how
 to find {\it sign-changing} solutions of some nonlinear elliptic equations,
 that is of importance in the  PDEs theory.   We would like to point
 out only these works in which the group and representation theory
 approach have been introduced and developed. This means that we do
 not mention many important papers and results based on a use of   this kind of  a symmetry approach.   Especially we do not expose  works  which use the same approach group theoretical approach,   but
are distinguished by their analytical form.

 Up to our knowledge the pioneer of this approach was  the paper \cite{BaWi} of Bartsch and Willem, where they used it  in  very particular
 form.
 They  studied a semilinear elliptic problem
 \begin{equation}\label{semilinear equation}
-\triangle u + b(\vert x\vert) \,u = f(\vert x \vert, u),\;\;  \
x\in \mathbb{R}^N, \;\;  u\in H^1(\mathbb R^N)\,.
 \end{equation}
The week solutions of \ref{semilinear equation} correspond to  the
critical points of functional
$$ \Phi(u):= \int_{\br^N} (\frac{1}{2} \vert \nabla (u) \vert^2 +
\frac{1}{2} b(\vert x\vert) \,u^2 - F(\vert x \vert, u)) \,dx \,,$$
with $F(r, u) = \int_0^u f(r,v) dv$ being the primitive of $f$ (cf.
\cite{BaWi}, also monographs \cite{KrRaVa}, \cite{Wi}). Observe that
the nonlinear functional $\Phi$ is $O(N)$-invariant with  respect to
the action of $O(N)$ on $\br^N$. A similar holds for the autonomous
nonlinear elliptic problem
\begin{equation}\label{autonomous problem}
-\triangle u  = f(u),\;\;  \ x\in  \Omega \subset \mathbb{R}^N, \;\;
u\in H^1(\Omega)\,
\end{equation}
where $\Omega $ is either bounded region with smooth boundary or
$\Omega = \br^N$ invariant with respect to the action of $O(N)$ on
$\br^N$. In this case the corresponding nonlinear functional is
equal to
$$ \Phi(u):= \int_{\Omega} (\frac{1}{2} \vert \nabla (u) \vert^2 - F(u)) \,dx \,$$
with $F(u) = \int_0^u f(v) dv$ being the primitive of $f$ (cf.
\cite{BaWi}, also monographs \cite{KrRaVa}, \cite{Wi}). $\Phi$ is
$O(N)$-invariant with  respect to the action of $O(N)$ on $\br^N$.

The existence of solutions of (\ref{semilinear equation}) and
(\ref{autonomous problem}) is obtained by standard and well-known
variational methods (cf. \cite{BaWi}, \cite{KrRaVa}, \cite{Wi}). In
particular, analytical assumptions  on  $f$, a  boundary condition,
and $\bz_2$-symmetry allow to apply the Ambrosetti-Rabinowitz
symmetry mountain pass theorem of \cite{AmRa} or  the fountain
theorem of \cite{Ba} (cf. \cite{KrRaVa},\cite{Wi} for general
references) which gives infinitely many solutions of discussed
problem (\ref{autonomous problem}). Note that any function space $E$
related to it has a natural $O(N)$ linear action given by
$$ g \, u(x):= u(g^{-1}x) \,.$$
Consequently, with every closed subgroup $G\subset O(N)$ we can
associate the corresponding linear subspace  $E^G$ of fixed points
of $G$, which is infinite-dimensional for the  discussed function
spaces.

Moreover, by posing  the problem (\ref{autonomous problem}) in
$E^G$, i.e restricting $\Phi$    to $E^G$,  finding its critical
points and finally using the  Palais  principle of symmetry (cf.
\cite{Palais}) we get solutions, or respectively infinitely many
solutions of studied problem which possess  the given symmetry. In
particular taking $G=SO(N)$ we get the radial solutions  (see
\cite{BaWi} for references up to $1993$).

In order to get non-radial solutions of (\ref{semilinear equation})
Bartsch and Willem took the subgroup $H \subset O(N)$ defined as
follows:
$$ \forall \;\; 2\leq m\leq  N/2, \;\;{\text{and}}\;\; 2m \neq N-1\;\; {\text{let}}\;\ H:= O(m)\times O(m)\times O(N-2m)\,.$$
 Next, they extended $H$ to $G
=\langle H\cup \{\tau\}\rangle$, where $\tau \in O(N)$ is the linear
map  transposing  the two first coordinates $ (x,y,z) \mapsto
(y,x,z)$. Observe that $\tau \in N(H)$, the normalizer of $H$ in
$O(N)$. This allows to define a representation $\rho$ of $G$ in a
function space $E$ putting \begin{equation}\label{interwinding} (g
f)(x) := \rho(\tau h) f(h^{-1}\tau^{-1}x) =  -
f(h^{-1}\tau^{-1}x),\; \; {\text{for}}\;\; g=\tau h \in G\,,
\end{equation}
where $\rho:G\to O(1)= \{-1,1\}$ is a representation of $G$ defined
by $\rho(\tau)=1$, $\rho(h)=1$ for every $h\in H$. Such functions,
thus solutions,  are called $\rho$-intertwining (cf. \cite{BrClMa}).

 It is easy to verify that
functions belonging to $E^G$ are sign changing with the zero set
containing  a hyperplane thus must not be radial. Consequently,
restricting $\Phi$ to $E^G$ they  have got infinitely many non
radial sign-changing solutions $\{u_n^m\}_{n=1}^\infty$  of
(\ref{semilinear equation}). We must add that they have to use the
space $E^G \subset E^H$ also for an analytical reason. Indeed, since
the embedding of $E^H\subset L^s(\br^N)$ is compact due the Lions
theorem, the restriction
 $\Phi_{|E^G}$ of functional $\Phi$  satisfies the Palais condition. Note that
the above trick  and assumption of the Lions theorem require $m\geq
2$, and $N- 2m\geq 2$ which can be satisfied only if  $N=4$, or
$N\geq 6$. Taking any $l\in I_N:=\{i\in \mathbb N:2\leq i\leq N/2,\,
2i\neq N-1\}$ as $m$ they get an infinite sequences $\{u_l^k\}$ of
solutions of (\ref{semilinear equation}) and showed that these
sequences of solutions are geometrically different with respect to
the symmetry, i.e. none of $u_n^m $ is in the $O(N)$-orbit of any
$u_n^k$ if $m\neq k$. In other words two solutions $u$, $v$ are
geometrically distinguished with respect to the action of group
$O(N)$ if
\begin{equation}
\sim \; \exists \; g\in O(N)\; {\text{such that}} \;\  v(x) =
u(g^{-1} \,x)\,.
\end{equation}
Consequently,  the number of those sequences of solutions containing
elements in different $O(N)$-orbits is at least
$\left[\log_2\frac{N+2}{3}\right],$ as shows a careful inspection of
\cite[Proposition 4.1, p. 457]{BaWi}.(here $[r]$ denote the integer
part of $r$).

The use of various  function spaces distinguished by their
symmetries property has been displayed in \cite{BaLiWe3} where the
authors made a general remark   that the amount of geometrically
distinct series of solutions obtained in this way is related  to the
number of partitions of $N$. They studied another nonlinear problem
which could be handled by this approach:
$$ \begin{cases} (-\triangle )^m \,u = \vert u\vert^{q-2}\, u, \;\;\;\; {\text{in}}\;\; \mathbb{R}^N\,,\cr
u \in \mathcal{D}^{m,2}(\mathbb{R}^N)\,
 \end{cases}
 $$
where $\mathcal{D}^{m,2}(\mathbb{R}^N)$ is a completion of
$C^\infty_0(\mathbb{R}^N)$ in a norm (cf. \cite{BaLiWe3}), $N>2m$,
and $q=\frac{2N}{N-2m}$.

The next step in a developing group theoretical scheme for finding
changing-sign solutions of nonlinear elliptic problems has been done
in \cite{KrMa}. For two groups of the form $H= O(N_1)\times
O(N_2)\times\cdots \times O(N_r) \subset O(N) $ and $K= O(M_1)\times
O(M_2)\times\cdots \times O(M_s) \subset O(N) $, $\sum_1^r N_i =
\sum_1^s M_j = N$,  there was given a  sufficient condition on
partitions $N_1, N_2, \dots, \, N_r$ and $M_1, M_2, \, \dots\, M_s$
which guarantees that the group $G=\langle H, K\rangle $ generated
by $H$ and $K$ acts transitively on the sphere $S(\br^N)$. It states
that there is not $\tilde{r} < r$ and $\tilde{s}< s$ such that
\begin{equation}\label{sufficient for transitivity}
\sum_{i=1}^{\tilde{r}} \, N_i = \sum _{j=1}^{\tilde{s}} \, M_j =
N^\prime < N\,.
\end{equation}

To assign to each such  subgroup $H$ as above a subspace  consisted
only of functions changing signs (if are nonzero!) we have to assume
that there exists  an element $\tau \in \mathcal{N}(H)\setminus H$,
where $\mathcal{N}(H)$ the normalizer of $H$ in $O(N)$, of order
two. The latter is satisfied if there exist $i \neq j $ such that
$N_i=N_j$. Consequently, every such  subgroup $H$ as above defines
an infinite series of changing-sign solutions in $E^G$ with
$G=\langle H, \tau \rangle$. Next, the condition (\ref{sufficient
for transitivity}) implies that for two subgroups $H$, $K$ the
corresponding subspaces $E^G$, $E^{G^\prime}$, with $G^\prime=
\langle K, \tau^\prime\rangle$, are linearly independent, because
$E^G\cap E^{G^\prime}$ consists of radial functions only, thus the
zero here.  In \cite{KrMa} we showed that there exists at least
$s_N={\rm card}\ I_N=\left[\frac{N-3}{2}\right]+(-1)^N$ different
pairs of subgroups satisfying  \ref{sufficient for transitivity}, by
an effective construction  of a special partition of $N$.  The
number $s_N$ does not depend only on the space dimension $N$ but on
the amount of constructed partitions.  Note also that $s_N\sim N/2$
as $N\to \infty,$ but the sequence $\{s_N\}_{N\geq 4}$ is not
increasing. This estimate is  affected by the fact that to apply the
Lions theorem on the compact embedding \cite{Lions} we had to assume
that for all $i,\, j$  $N_i, \,M_j \geq 2$ in a constructed
partition associated with $H$, respectively $K$.  This way we proved
the existence of $s_N$ sequences of non-radial, sign-changing weak
solutions  such that elements in different sequences are mutually
distinguished by their symmetry properties. In \cite{KrMa} we
considered particular problem
$$\left\{
\begin{array}{lcl}
-\triangle_p u +|u|^{p-2}u=  \mathcal{K}(x)f(u),& \ x\in \mathbb{R}^N,\\
u\in W^{1,p}(\mathbb R^N),
\end{array}\right.\eqno{{\rm (P)}}$$
 when $p>N$, the space dimension $N$ is large enough,  and $f$ has an
oscillatory behavior at the origin. Here,  $\triangle_p u={\rm
div}(|\nabla u|^{p-2}\nabla u)$ is the usual $p$-Laplacian of $u$,
 $\mathcal{K}:\mathbb R^N\to \mathbb R$ is a  measurable function, and
$f:\mathbb R \to \mathbb R$ is continuous.

This group theoretical scheme can be applied to study similar
problems with different analytical part which have been successfully
used by A. Kristaly and co-authors in a couple of papers (see
\cite{KrRaVa} for references).

However, there are still remained open questions and objects for
studies: \vskip 0.5 cm
\newpage

\begin{enumerate}\label{questions}

\centerline{ {\bf Q \ref{questions}.}  {\bf Questions}}

\item{Analyze whether it is possible to  improve the construction
given in \cite{KrMa} to enlarge the number of subgroups for which
the spaces  $ E^H$  of fixed points in $E$ can be used for the
construction of subspaces described above. In particular get rid of
a fixed order of partitions $(N_1,N_2,\, \dots\,,N_r)$ and $(M_1,
M_2,  \, \dots\,,M_s)$ in condition (\ref{sufficient for
transitivity}).}

\item{Examine a dependence of a choice of such a subgroup on the orthogonal splitting into subspaces  of $\br^N$. In particular
what happens  if  one of the groups $H,\, K$ as above is constructed
with respect  to  a partition  of $N$ which is taken  the canonical
orthonormal basis but the second with respect to a partition
corresponding to another orthonormal  basis.}
\item{Find a largest number $s_N$ of subgroups $H \subset O(N) $ with
the  normalizers  $\mathcal{N}(H)$  containing  involutions  in
$\mathcal{N}(H)\setminus H$ such that for $H\neq K$ and $G=\langle
H,\tau\rangle $ and $G^\prime =\langle K,\tau^\prime \rangle $ we
have
 $$ E^G \cap E^{G^\prime} =\{0\}, \;\;{\text{and for the orbit}}\;{\text{of}}\; 0\neq v \in E^G\;\;{\text{we have}}\;\;   O(N)\,
 v \cap  E^{G^\prime} =\emptyset\,.$$}
 \item{Get an information about the nodal set of every function
 $u\in E^G$.}

\end{enumerate}

In this work we give answers to all the  questions formulated in {Q
\ref{questions}}. More precisely, the answer to Question (1) is
positive and is contained in considerations of  Sections 2, 3, and 4
(e.g. Theorem \ref{equivalent orthogonal flags}, Proposition
\ref{adaptation of KrMa}). An answer to Question (2). states that  a
choice of orthogonal basis determining each of these two groups does
not  have an affect on  this property. Only a relation between
 the partitions $(N_1,N_2,\, \dots\,,N_r)$ and
$(M_1, M_2, \, \dots\,,M_s)$ is essential (cf. Theorem \ref{main}).
Next, by reducing the problem of  estimating the amount $s_N$ of
pairs of subgroups $(H, K)$ with this property to an estimate of
number of not equivalent special partitions of $N$ we show that the
rate of growth of $s_N$ is exponential (cf. Theorem \ref{final
theorem}) which is an answer to Question (3). Finally, Corollary
\ref{nodal set} gives an answer to Question (4).

{\begin{rema} Since we studying ,  $O(N)$ invariant problems if a
function $ u $ is a solution then $(gu)(x) = u(g^{-1}x)$ is also a
solution. Consequently, we are interested in solutions (series of
solutions) which are geometrically distinct. The latter means that
none of the solutions from one series is in the $O(N)$ orbit of a
solution from another series.
\end{rema}}

 \vskip 0.5cm

The paper is organized as follows. After the introduction and  short
opening on the fixed point spaces of a representation, in Section
\ref{subspaces} we provide an information on the groups  acting
transitively on spheres called the Borel groups. Next we  introduce
the notion of orthogonal  Borel subgroups, or correspondingly the
maximal orthogonal subgroups (cf. Definitions \ref{Borel subgroup},
\ref{Borel subgroup}, \ref{maximalBorel subgroups}) and relate them
to partial orthogonal flags in $\mathbb{R}^N$ in Proposition
\ref{form of Borel subgroup}. In Section \ref{equivalent maximal
orthogonal Borel subgroups} we discuss the action of $O(N)$ on the
set of all partial orthogonal flags and the set of all maximal
orthogonal Borel subgroups. We show that  these actions coincide
which implies the correspondence of the orbits (Proposition
\ref{equivalents of actions}. This gives a combinatorial description
of these equivalence classes as the partition of the number $N$
(Theorems \ref{equivalent orthogonal flags}, \ref{combinatorially
equivalent partitions}). Next we derive the Weyl group of a partial
flag (Proposition \ref{Weyl group of flag}), thus the Weyl group  of
a Borel subgroup (Corollary \ref{Weyl group of subgroup}). We end
this section with a survey of information on the amount of
partitions of given  number $N$, next its partitions without a
repetition, and finally the amount of partitions with every summand
$\geq 2$. We  conclude with Proposition \ref{special partitions} in
which we effectively construct $s_N$ partitions with nontrivial Weyl
groups and each summand $\geq 2$ in the amount $s_N$ growing
exponentially in $N$. In the last section we prove our main Theorem
\ref{final theorem} which proof is based on  Proposition \ref{linear
independence of fixed points}. The latter states that each  class of
$s_N$ constructed partitions of $N$ determines a subspace of the
functional space $E$ in such a way that it is the fixed points space
of a subgroup $\tilde{H} \subset O(N)$,  and subspaces corresponding
to different partitions are linearly independent, i.e. their
intersection is equal to $\{0\}$. Restricting the functional to each
of these subspaces we get $s_N$ infinite series of geometrically
distinct solutions.  The thesis of Theorem \ref{final theorem}
automatically applies to the  results of \cite{BaWi},
\cite{BaLiWe3}, \cite{KrMa}, \cite{Kr}, and all  similar quoted in
\cite{KrRaVa}, giving a much  larger number of series of
geometrically distinct solutions than  of those  papers. Finally, in
Section \ref{Appendix} we provide additional information about the
types of partitions we use and show describe subgroups generated by
two maximal orthogonal Borel subgroups in $\mathbb{R}^N$ associated
with two partitions of $N$ (Proposition \ref{KrMrTheorem}, Theorem
\ref{group generated by aprtitions}).

  \zz{In Section
\ref{fotetel-J} we will formulate our main result. In Section
\ref{variational-sect} we
  prove the existence of a sequence  of  weak solutions of (P) for a
suitably chosen subspace of $W^{1,p}(\mathbb  R^N)$ which are
sign-changing
  non-radial by the symmetry property of the subspace.
 In the last section, we construct $s_N$ special
 subspaces of $W^{1,p}(\mathbb R^N)$ for which one can apply the arguments
 from Section  and these subspaces have only the 0
 as a common element. In this way, we produce $s_N$  sequences of solutions
 of (P) which  do no contain similar elements from symmetrical
 point of view.}

\section{Subspaces of fixed points of subgroups and Borel subgroups}\label{subspaces}
 By $O(N)$ denote the group of all linear orthogonal maps of the
Euclidean space of dimension $N$ and by  $SO(N)\subset O(N)$ its
connected component of $\{e\}$ consisting of maps preserving
orientation.

Let $\Omega \subset \br^N$ be an open region in $\br^N$ which, is
$O(N)$ invariant, i.e. $\Omega = D_r^N$ if it is bounded, or $\Omega
= \br^N$ if it is unbounded.  From now on, by a functional space $E$
with a domain $\Omega $ we understand a completion of
$C^\infty(\Omega)$ or $C^{\infty}_0(\Omega)$ in any linear topology,
e.g. any topology induced by a norm, such that the natural linear
action

$$  O(N) \times E \to E, \;\; (g,u)\mapsto u(g^{-1}x)$$
is continuous and preserving norm if it is the case.

 Let $H, \,K
\subset O(N)$ be two closed subgroups and $E^H$, respectively $E^K$
the fixed point subspaces.

By a general theory of transformations of groups we have
$$ u\in E^H \; \Longleftrightarrow gu(x)= u(g^{-1} x) \in
E^{gHg^{-1} } \,.$$ In particular, if $K= gHg^{-1}$ then every
element $u\in E^H$ has the same orbit as  $v=gu \in E^K$.

Let $\mathcal{G}(N)$ be the set of all closed subgroups of $O(N)$
with the action of $O(N)$ by conjugation, and next  $\mathcal{S}(N)$
be the set of all conjugacy classes of  (closed!) subgroups  of
$O(N)$ and $\{H_s\}_{s\in \mathcal{S}}$ a complete set of
representatives of $\mathcal{S}$.

 Since we are
interested in finding geometrically distinct solutions there is
reasonable to take only one representative  $H_s \in \mathcal{G}(N)$
of the class  $ [H_s] \in \mathcal{S}(N)$.

\begin{rema}
We must emphasize that $[H_s] \neq [H_{s^\prime}]$  $ \nRightarrow $
$E^{H_s}\cap E^{H_{s^\prime}} = \{0\}$ in general.
\end{rema}

\zz{In next we will show that for a  special family of  subgroups,
called the Borel subgroups, the intersection consists of functions
radial on a subspace $V\subset \br^N$, and if the Weyl groups   of
them are non-trivial and preserve $V$ then the intersection is equal
to $\{0\}$, i.e. these two subspaces are linearly independent.}

\zz{\section{Borel subgroups}\label{section Borel subgroups}}

Next, we will define  a class of subgroups of $O(N)$,  respectively
$SO(N)$, called the Borel subgroups. Beforehand, we present a short
opening to justify the  name this  notion.
\begin{defi}\label{Borel group definition}  We call a closed  subgroup $G\subset O(N)$ an {\underline{
orthogonal Borel group}}, if  $G$ acts transitively on the unit
sphere $S(\mathbb{R}^N)= S^{N-1}$.  If $G$ is connected, i.e. $G
\subset SO(N)$    then we call $G$ {\underline{connected orthogonal
Borel group}}.

\end{defi}

Note that $G$ is an orthogonal Borel group iff it acts transitively
on the sphere of any radius $r>0$.

Also we have to point out that the above notion of the
{\underline{orthogonal Borel group}} is different that a notion of
the {\underline{Borel group}} used in the algebraic geometry. An
expiation of the origin of this notion and an information on the
corresponding Borel theorem is included to the Section
\ref{Appendix}.
\begin{rema}{
If a group $G$ acts transitively on a $G$-space $X$ then $X$ is
($G$-equivariantly) homeomorphic to the orbit $Gx$ of any point
$x\in X$, i.e. it is homeomorphic to the homogenous space $G/H$,
with $H=G_x$. If $G$ acts smoothly, or by isometries respectively
then $X$ is   diffeomorphic,  or  correspondingly isometric to
$G/G_x$. }
\end{rema}

 \begin{rema} Our assumption  of the Definition of  orthogonal Borel group, that $G\subset O(N)$,  automatically guarantees that it acts effectively.
\end{rema}

\begin{rema}\label{action of component of identity} Note that if $G_0\subset G$ acts transitively on
$S(\br^N)$ then $G$ does too. For $ N\geq 2$, if $G$ acts
transitively on $S(\br^N)$ then $G_0$ does too. Indeed, otherwise
$S(\br^N)$ would be a finite ($\vert G/G_0 \vert$) union of closed
submanifolds which is impossible.

This means that that for $N\ge 2 $ a group $G\subset O(N)$ is a
Borel group if and only if its connected component $G_0\subset
SO(N)$ is a connected Borel group.
\end{rema}

Let $H \subset O(N)$ be a closed subgroup. Then the  linear space
$\mathbb{R}^N$ is the space  of  orthogonal representation of $H$
given by the restriction to $H$ the canonical representation of
$O(N)$ in $\mathbb{R}^N$.

We denote by $(V, \rho_H)$, or shortly $V$, the space $\mathbb{R}^N$
with the  defined  above representation structure $\rho_H: H\to
O(N)$.

\begin{defi}\label{Borel subgroup}
We call a closed   subgroup $H\subset O(N)$ an
{\underline{orthogonal}} {\underline{Borel subgroup}} if:

\begin{itemize}

\item[$A_1)$] {The restricted representation $(V,H)$ of $H$ in
$\mathbb{R}^N$ decomposes into a direct sum
$V={\underset{j=1}{\overset{r} \oplus}}\, V_j$, $r\geq 2$,  of
pairwise non-isomorphic irreducible representations $V_j$  of $H$.}

\item[$A_2)$] {For every $1\leq j \leq r$ the group $H$ acts
transitively on $S(V_j)$.}

\end{itemize}
\end{defi}
\begin{defi}\label{connected Borel subgroup}
We call a closed  connected  subgroup $H\subset O(N)$ a
{\underline{connected orthogonal}} {\underline{Borel subgroup}} if
it is an orthogonal Borel subgroup and is connected.
\end{defi}
Note that  for any connected Borel subgroup we have $H\subset SO(N)
$.
\begin{rema}\label{dimension of conneceted Borel subroups}\rm Since $O(1)=\bz_2=\{-1,1\}$, but  its connected
component $SO(1)$ does not act transitively on $S(\br)=\{-1,1\}
\subset \br $,  for all irreducible representations  $V_j$ of
Definition \ref{connected Borel subgroup}  we have to assume that
\begin{equation}\label{dimension assumption}
 N_j =
\dim_\br V_j \geq 2\,.
\end{equation}

\end{rema}

Let $H$ be an orthogonal Borel subgroup as in Definition \ref{Borel
subgroup} and $N_j=\dim V_j$ the dimension of $V_j$. Note that
${\underset{j=1}{\overset{r} \sum}}\, N_j= N$. Furthermore, for
every $1\leq j\leq r$ we have canonical projection $p_j : V\to V_j$
given by
$$ p_j(v):= \int_H \rho_H(h) \, \chi^*_j(h) d\mu(h)\,,$$
where $\mu$ is the Haar measure on $H$,  $\chi_j: H \to \mathbb{R}$
is the character of irreducible representation $V_j$, and $\chi^*$
denote the character of representation conjugated to $\chi$.
Moreover $V_j=p_j(V)$, and we have  a canonical embedding
$\epsilon_j: V_j\hookrightarrow V$, $v \mapsto (0, \,0, \, \cdots,
\,v,\,0,\, \cdots,\, 0)$.

 Since $H \subset O(N)$, every  $p_j$ is an orthogonal projection, $\epsilon_j$  is an orthogonal embedding,  and they are
 $H$-equivariant.

 Consequently for any $h\in H$ the formula $ h \mapsto \rho_j(h)$
defines a homomorphism, denoted also  $p_j$, from $H$ into $O(N_j)$
which can be defined as

\begin{equation}\label{projections of H}
 p_j(h) (v):=  p_j(\,\rho_H(h)(\epsilon(v))\,)\,, \;\; {\text{where}}\;\;v\in
V_j\simeq \mathbb{R}^{N_j}\,.
\end{equation}

 \vspace{-0.0cm}

 Denote by $H_j= p_j(H)$ the image of $H_j$ in $O(N_j)$.

 \begin{lema}\label{projection of Borel subgroup}

 For every $1\leq j\leq r$ the group $H_j=p_j(H)$ is an orthogonal Borel group
 of $O(N_j)$.
 \end{lema}

\noindent {\sc{Proof.}} By definition  $H_j \subset O(N_j)$.
 On the other hand, by assumption $A_2)$ $H$ acts transitively on $S(V_j)$,
 but the action of $H$ on $S(V_j)$ factorizes through $H_j$,
 i.e.
 for every $v_j=(0, \cdots, \, v,\,\cdots, 0)\in S(V_j)$ we have
 $\rho_j(h)(v_j)= \rho_j(p_j(h))(v_j)$. Indeed $ \rho_H(h)=
 p_j(\rho_H(h)(v_j))$, since $\rho_H$ preserves each $V_j$, and
$p_j(\rho_H(h)(v_j)) =  p_j\rho_H(h)(\epsilon(v)) = p_j(h)(v)$. \hbx

This leads to the following characterization of the orthogonal Borel
subgroups of $O(N)$.

\begin{propo}\label{form of Borel subgroup}
Every orthogonal  Borel subgroup $H$ of $O(N)$ is isomorphic to a
product $H \simeq H_1\times H_2 \times \, \cdots\, \times H_r$ where
$H_j$ is an orthogonal  Borel group of $O(N_j)$ and
${\underset{j=1}{\overset{r} \sum}}\, N_j= N$.

More precisely, if the restricted representation $(V,\rho_H)$, $\dim
V=N$, of $H$ decomposes into $r$ irreducible summands
$V={\underset{j=1}{\overset{r} \oplus}} V_j$, with $\dim V_j =N_j$,
then $H_j=p_j(H)$ (see {\rm{(\ref{projections of H})}}).
\end{propo}

\noindent {\sc{Proof.}} The product map $p:= p_1 \times\,
\cdots\,\times p_r: \, H \to H_1\times \, \cdots\,\times H_r$ is a
homomorphism as the product of homomorphisms, and is onto by the
definition of $H_j$.

It is enough to show that it is injective. Indeed, let $p(h) = e
=(e,\, \cdots\,, e) \in H_1\times \, \cdots\,\times H_r$, i.e.
$(p_1(h),\,\cdots\, p_r(h))(v_1, \,\dots\,,v_r)=(v_1, \,\dots\,,v_r)
$ for every $v=(v_1, \,\dots\,,v_r)\in {\underset{j=1}{\overset{r}
\oplus}} V_j$. But $V={\underset{j=1}{\overset{r} \oplus}} V_j$, and
$\rho_H(h)(v)= {\underset{j=1}{\overset{r} \sum}}
\rho_j(p_j(h))(v_j)$ as we have shown in the proof of Lemma
\ref{projection of Borel subgroup}. This shows that $\rho_H(h)= {\rm
id}$ and consequently $h=e$, since $\rho_H$ is given by the
embedding of the subgroup $H$ in $O(N)$. \hbx
\begin{rema}
Observe that the assumption that $\rho_H$ is a direct sum of
pairwise non-isomorphic irreducible representations is necessary for
  the transitivity of action of $H$ on the spheres $S(V_j)$ of the
summands of this decomposition. Indeed, the canonical projections
$p_j$ (thus canonical decomposition of $V$) are projections onto the
isogenic subrepresentations of $V$, i.e. onto subrepresentations
which are multiplicities of a given irreducible representation
$\rho_j: H \to O(N_j)$. If $(V,\rho_H)$  contains a representation
$k V_j$, $V_j$ irreducible, $k>1$, then $H$ does not act
transitively on $S(kV_j)$, because it acts diagonally on $k
V_j={\underset{j=1}{\overset{k} \oplus}} V_j$.

\end{rema}

For a given orthogonal Borel subgroup $H\subset O(N)$ and the
restricted representation $(V,\rho_H)$,
$\;V={\underset{j=1}{\overset{r} \oplus}} V_j$ we put
$$ V^j:= {\underset{l=1}{\overset{j} \oplus}} V_l\,.$$
Note that $V_1=V^1$, $V^r=V$, $V^i \subset V^j$ if $i<j$, i.e. the
family $\{V^j\}_{j=1}^r$   forms a filtration of $V$ by  linear
subspaces. Put $n_j= {\underset{l=1}{\overset{j} \sum}} N_l$ equals
to the dimension of $V^j$. Additionally put $V^0=V_0 =\{0\}$.

Note that a filtration of $V=\br^N$ by $$\{0\}= V^0 \varsubsetneq
V^1 \varsubsetneq V^2\, \dots \varsubsetneq \, V^r=V=\br^N$$ is
called {\it a partial flag} of $V$. Furthermore, if $H\subset O(N)$
preserves this filtration $\{V^i\}$  (the flag) it preserves also
the subspaces of orthogonal complements $V_i= V^{i-1}_\perp  $ in
$V^i$. This leads to an equivalent definition of an orthogonal Borel
subgroup of $O(N)$.

\begin{prop}\label{another definition of Borel subgroup}
Let $H\subset O(N)$ be a subgroup. Then $H$ is an orthogonal  Borel
subgroup if and only if there exists a partial flag $V^1
\varsubsetneq V^2\, \dots \varsubsetneq \, V^r$ of $\br^N$ which is
preserved by $H$ and for every $1\leq i\leq r$ the induced action of
$H$ on the projective space $\mathbb{P}(V^i/V^{i-1})$ of quotient
space $V^i/V^{i-1}$ is transitive. \hbx
\end{prop}

The above justifies a use of the following notion. For every partial
flag  $$\{0\}= V^0 \varsubsetneq V^1 \varsubsetneq V^2\, \dots
\varsubsetneq \, V^r=V\;\;{\text{in}} \;\; V=\br^N $$ we call the
corresponding orthogonal decomposition
$$ V={\underset{j=1}{\overset{r} \oplus}} V_j, \; \;
{\text{where}}\;\; V_j\supset V^j = V^{j-1}_\perp $$ {\it the
orthogonal partial flag}.

\begin{rema}
Note that the isomorphism of $H \simeq H_1 \times H_2 \times
\cdots\times H_r$ of Proposition \ref{another definition of Borel
subgroup} depends on the flag $\{0\}= V^0 \varsubsetneq V^1
\varsubsetneq V^2\, \dots \varsubsetneq \, V^r=V=\br^N$.
\end{rema}

\medskip

For our analytical considerations  we will use only a special class
 of orthogonal Borel subgroups, each of them is canonically associated with a
 partial flag.

 \begin{defi}\label{maximalBorel subgroups}
We will call $O(V_1) \times\, O(V_2)\, \times\, \dots \,\times
O(V_r) \equiv O(N_1) \times\, O(N_2)\, \times\, \dots \,\times
O(N_r)$ the {\underline{maximal orthogonal Borel subgroup}}, and
respectively $SO(V_1) \times SO(V_2) \times\, \dots \,\times
SO(V_r)$ the {\underline{maximal connected orthogonal  Borel
subgroup}}  associated with given partial orthogonal flag $\{V_1, \,
V_2,\,\dots\,, V_r\}$.
 \end{defi}

By $\mathcal{G}_B$ we denote the set of all maximal orthogonal Borel
subgroups, and respectively by $\mathcal{G}_B^0$ the set of all
maximal connected  orthogonal Borel subgroups.

\section{Equivalent maximal orthogonal Borel
subgroups}\label{equivalent maximal  orthogonal Borel subgroups}

Now we would like to describe a number of non-equivalent orthogonal
Borel subgroups that correspond to an action of the group $O(N)$, or
$SO(N)$ on the set of orthogonal partial flags. To shorten notation,
in  this section by an orthogonal, correspondingly connected
orthogonal, Borel group we mean  the maximal orthogonal Borel,
respectively maximal connected Borel subgroups .

 Let us denote by
$\mathfrak{F}_r(N)$ the set of  flags of length $1 \leq r \leq N$ in
$\br^N$.

\subsection{The action of $O(N)$ on the set of partial flags}

Observe that if $g\in O(N)$, in particular if $g\in SO(N)$, and
$V={\underset{j=1}{\overset{r} \oplus}} V_j$ is an orthogonal
partial flag then the family of spaces $ \{g \,V_j\}_{j=1}^r$ forms
also an orthogonal partial flag of the same length.
\begin{defi}\label{action on flags}\rm
The mapping  $(g,V={\underset{j=1}{\overset{r} \oplus}} V_j) \mapsto
\{g \,V_j\}_{j=1}^r$  defines an action of $O(N)$ or $SO(N)$ on $
\mathfrak{F}_r(N)$.
\end{defi}

The isotropy group of any such an orthogonal partial flag is equal
to $$O(V_1) \times\, O(V_2)\, \times\, \dots \, \times \,O(V_r)
\equiv O(N_1) \times\, O(N_2)\, \times\, \dots \, \times O(N_r)\,,$$
or respectively
$$SO(V_1) \times\, SO(V_2)\, \times\, \dots \,\times \, SO(V_r) \equiv SO(N_1)
\times\, SO(N_2) \times\, \dots \,\times SO(N_r)$$ if we consider
the action of $SO(N)$.

Consequently, the orbit of any orthogonal partial flag is equal to
$$O(N)/(O(N_1) \times\, O(N_2) \times\, \dots \,\times O(N_r)) \; {\text{or}}\; SO(N)/(SO(N_1) \times SO(N_2) \times\,
\dots \,\times SO(N_r))$$ or respectively.

Note that by Definition \ref{Borel subgroup} an orthogonal Borel
subgroup $H$ associated with given partial orthogonal flag
$V={\underset{j=1}{\overset{r} \oplus}} V_j$ is  contained in the
isotropy group of $\{V_j\}$, i.e. in $O(V_1) \times\, O(V_2)\,
\times\, \dots \,\times  O(V_r) \equiv O(N_1) \times\, O(N_2)\,
\times\, \dots \,\times  O(N_r)$.  By the definition, if it is a
connected orthogonal Borel subgroup it is a subgroup of $SO(N)$ thus
a subgroup of the isotropy group $SO(V_1) \times SO(V_2) \times\,
\dots \,\times SO(V_r) \equiv  SO(N_1) \times SO(N_2) \times\, \dots
\,\times SO(N_r) $.

 On the other hand we have the   orthogonal Borel subgroup, or correspondingly the   connected orthogonal Borel subgroup,  of
 the partial orthogonal flag $\{g V_j\}_{j=1}^r$ is equal to  $O(gV_1) \times\, O(gV_2)\, \times\, \dots \, \times O(gV_r))= g
O(V_1)g^{-1} \times g O(V_2)g^{-1} \times\, \dots \,\times g
O(V_r)g^{-1}\equiv O(N_1) \times O(N_2) \times\, \dots \,\times
O(N_r) $
 or respectively
  $SO(gV_1) \times SO(g V_2) \times\, \dots \,\times SO(g V_r))= g
SO(V_1)g^{-1} \times g SO(V_2)g^{-1} \times\, \dots \,\times g
SO(V_r)g^{-1}\equiv SO(N_1) \times SO(N_2) \times\, \dots \,\times
SO(N_r) $.

Summing up, we have the following proposition.
\begin{propo}\label{equivalents of actions}
 The natural
action of the group $O(N)$, correspondingly $SO(N)$, on the set
$\mathfrak{F}_r(N)$ of all orthogonal partial flags of  length $r$
given by  $( g, \{ V_j\}_{j=1}^r) \mapsto \{g V_j\}_{j=1}^r$
corresponds to the action by conjugation $(g, H) \mapsto gHg^{-1}$
on the set $\mathcal{G}_B(N)$, correspondingly  $\mathcal{G}_B^0(N)$
of   orthogonal Borel subgroups, respectively the   connected
orthogonal Borel subgroups, associated with them.

\end{propo}

\begin{defi}\label{equivalence of partial flags}
\rm We say that two orthogonal partial flags $\{ V_j\}_{j=1}^r$ and
$\{ W_i\}_{i=1}^s$ are equivalent if they are in one orbit of the
above action, i.e. if there exists $g\in O(N)$, respectively  $g\in
SO(N)$, such that $$g \,\{ V_j\}_{j=1}^r = \{ W_i\}_{i=1}^s\,.$$
\end{defi}

Now we show a condition which is necessary and sufficient for the
equivalence of orthogonal partial flags. To do it we need new
notation. We say that a sequence $N_1, N_2, \, \dots\,N_r$ of
natural numbers such that $N_1+N_2+\,\dots \,+N_r=N$ is a {\it
partition}  of $N$ and $r$ is the {\it length} of this partition.

\begin{defi}\label{definition equivalent partitions} \rm
We say the  two partitions $N_1, N_2, \, \dots\, N_r$ and
 $M_1, M_2, \, \dots\, M_s$ of $N$ {\underline{are equivalent}}
if $r=s$ and there exists  a permutation  $\sigma$ of the set $\{1,
2\, \dots\,, r\}$ such that $N_j=M_{\sigma(j)}$.
\end{defi}
Note that  any partition $N_1, N_2, \,\cdots \,, N_r$ of $N$ is
equivalent to a partition $\bar{N}_1, \bar{N}_2,\,  \cdots \, ,
\bar{N}_r $  such that $\bar{N}_1\leq \bar{N}_2\leq \, \cdots\, \leq
\bar{N}_r$ , in particular if $N_j \neq N_i$ for $i\neq j$ then the
partition $N_1, N_2, \, \dots\, N_r$ is equivalent  to a partition
such that  $ \bar{N}_1 < \bar{N}_2 < \,\cdots\, < \bar{N}_r$.

We are in position to formulate the main result of this section.
\begin{teo}\label{equivalent orthogonal flags}
 Two orthogonal partial flags $\{ V_j\}_{j=1}^r$
and $\{ W_i\}_{i=1}^s$ are equivalent if and only if $r=s$ and  the
corresponding to them partitions of $N$:  $(N_1, N_2, \, \dots\,
N_r)$ and  $(M_1, M_2, \, \dots\, M_s)$,  with $N_j=\dim V_j$,
respectively $M_i=\dim W_i$, are equivalent.
\end{teo}

   \noindent{\sc Proof.}
Since an orthogonal linear map $g$ preserves the dimension of spaces
and their orthogonality the partitions of $N$ associated with $\{
V_j\}_{j=1}^r$ and with $ \{g\, V_j\}_{j=1}^r$ are the same which
proves the necessity of condition.

 Now suppose that $r=s$ and there exists a permutation $\sigma$ of
the set $\{1, 2\, \dots\,, r\}$ such that $N_j=M_{\sigma(j)}$.
Changing indices at the partition $(M_1, M_2, \, \dots\, M_s)$ and
corresponding orthogonal partial flag $\{ W_i\}_{i=1}^s$ we can
assume that $N_j=M_j$ for $1\leq j\leq r$. For $1\leq j \leq r$ let
$\{\vec{e}_i^j\}$, $i=1, \,\dots\,,N_j$ be an orthonormal basis of
$V_j$, and $\{\vec{f}_i^j\}$, $i=1, \,\dots\,,N_j$ the analogous
basis of $W_j$. From linear algebra, it follows that there exists an
orthogonal map $g \in O(N)$ such that $$ g(\vec{e}_i^j)
=\vec{f}_i^j\,.$$ Obviously, $g (V_j)\subset W_j$, thus $g V_j =W_j$
for every $ 1\leq j\leq r$.  Moreover, fixing orientations in each
$V_j$ and $W_j$ we can choose each basis $\{\vec{e}_i^j\}$ and
$\{\vec{f}_i^j\}$ consistent with their orientation. Then, $g\in
SO(N)$ and $g_{V_j}: V_j \to W_j$ is also a preserving orientation
linear orthogonal map. This shows that this condition is sufficient
for the equivalence of orthogonal partial flags. \hbx

As we stated in Section \ref{subspaces} we are interested in a
choice of one   orthogonal  Borel subgroup from each equivalence
(conjugacy) class. By Proposition \ref{equivalents of actions} and
Theorem \ref{equivalent orthogonal flags} we get the following
statement
\begin{coro}\label{final equivalence}
Two    orthogonal Borel subgroups $H=O(V_1) \times\, O(V_2)\,
\times\, \dots \,\times O(V_r) \equiv O(N_1) \times\, O(N_2)\,
\times\, \dots \,\times O(N_r)$ and $H^\prime =O(V^\prime_1)
\times\, O(V^\prime_2)\, \times\, \dots \,\times O(V^\prime_s)
\equiv O(N^\prime_1) \times\, O(N^\prime_2)\, \times\, \dots
\,\times O(N^\prime_s) $   are equivalent if and only if $r=s$ and
the partitions  of $N:$  $N_1, N_2, \,\dots\,,N_r$, and $N^\prime_1,
N^\prime_2, \,\dots\,,N^\prime_r$ are equivalent.

\end{coro}
In other words, to chose one representative of each equivalence
class of the   orthogonal Borel subgroups  is enough  to fix an
orthogonal basis $\mathcal{E}=\{e_1, \,e_2, \,\dots\,,e_N\}$, next
fix a representative $N_1,N_2, \,\dots\,, N_r$  of each equivalence
class of partitions, and finally construct a partial orthogonal flag
$\{0\} \subset V^1 \subset V^2\subset \cdots\,\subset V^r=\br^N$,
$V^i={\rm {span}} \{e_1,e_2, \,\dots\,,e_{\bar{N}_j}\} $, with
 $\bar{N}_j=
\sum_1^j N_j$,  defining  the   orthogonal subgroup
$$H=O(V_1) \times\, O(V_2)\, \times\, \dots \,\times O(V_r) =
O(N_1) \times\, O(N_2)\, \times\, \dots \,\times O(N_r)\,.$$ \hbx

\subsection{A combinatorial description of the  action}
Now we would like to give a combinatorial condition for two
partitions of $N$ to be equivalent in the sense  of Definition
\ref{definition equivalent partitions}

Let $\pi_r(N)=\{N_1, N_2, \, \dots\, N_r\}$ be a partition of $N$,
denoted shortly  by $\pi$. The  set of all partitions of $N$ of
length $r$  we denote by $\Pi_r(N)$.

 The permutation group of $r$-symbols $\mathfrak{S}(r)$
acts on $\Pi_r(N)$ by permuting the indices. Observe that two
partitions $\pi, \, \pi^\prime \in \Pi_r(N)$ are in one orbit of the
action of $\mathfrak{S}(r)$ iff they are equivalent in the sense of
Definition \ref{definition equivalent partitions}. On the other hand
they are in the same orbit of the  action of $\mathfrak{S}(r)$ on
$\Pi(r)$ iff the isotropy groups $\mathfrak{S}_\pi$ and
$\mathfrak{S}_{\pi^\prime}$ are conjugated in $\mathfrak{S}(r)$  as
follows from general theory of actions of groups.

Note that for a partition $\pi=\{N_1, \, \dots \,N_r\}$ with
$N_i\neq N_j$ for $i\neq j$ we have $\mathfrak{S}_\pi = e$ the
identity permutation.

Our task is to describe the  isotropy group $\mathfrak{S}_\pi$ of a
partition $ \pi \in \Pi_r(N) $. To do it we define a function
\begin{equation}\label{definition phi}
    \phi_\pi: \{1, \, 2, \,\dots\, , \, N\} \to 2^{\{0, \,1,\, \dots\,
    r\}}
\;\; {\text{defined as}} \;\;
    \phi_\pi(n) = \{1\leq j \leq r:\, N_j = n \}
\end{equation}
Note that $\phi(k) \cap \phi(n) =\emptyset$ if $k\neq n$.

Now we  define next function
\begin{equation}\label{definition psi}
    \psi_\pi: \{1, \, 2, \,\dots\, , \, N\} \to {\{0, \,1,\, \dots\,
    r\}}
    \;\; {\text{defined as}} \;\;
    \psi_\pi(n) = \vert \phi_\pi(n) \vert
\end{equation}
where $\vert A \vert $ denotes the cardinality of a finite set $A$
with a convention that $\vert \emptyset \vert =0$.

\begin{lema}\label{isotropy of partition}
Let $\pi=\{N_1, \, \dots \,N_r\}$ be a partition of $N$ such that
there is  $q$   different values $n_1,\, \dots\, n_q$ in the
sequence $N_1, \, N_2,\, \dots\,, N_r$. Then the isotropy group
$\mathfrak{S}_\pi \in \mathfrak{S}(r)$ of the action of
$\mathfrak{S}(r)$ on $\Pi_r(N)$ is equal to $
\mathfrak{S}(\phi_\pi(n_1)) \, \mathfrak{S}(\phi_\pi(n_2)) \,
\cdots\, \mathfrak{S}(\phi_\pi(n_q))\,.$
\end{lema}\label{isotropy gropu of permutation}
Note that since the supports of permutations in
$\mathfrak{S}(\phi(n_i))$, with $\phi_\pi(n_i)\neq 0$, are disjoint
for different $i$, all they commute. Consequently,
\begin{equation}
     \mathfrak{S}_\pi =  \mathfrak{S}(\phi_\pi(n_1)) \times \mathfrak{S}(\phi_\pi(n_2)) \times
    \cdots\, \times
\mathfrak{S}(\phi_\pi(n_q)) \;\; \text{is a group of order} \;\;
 \psi_\pi(1) ! \, \psi_\pi(2) !
    \cdots\,
\psi_\pi(N)!
\end{equation}
with the  usual convention that $0!=1$.

Note the length $r$ of any abstract  partition should be $\leq N$,
and under our assumption that $\dim V_j=N_j\geq 2$ we have
 $r\leq \frac{N}{2}$.

Now we are in position to formulate a theorem which characterizes
combinatorially  equivalent partitions of $N$ of the same length
$r$.

\begin{teo}\label{combinatorially equivalent partitions}
Two partitions $\,\pi =\{N_1, \, \dots \,N_r\} $ and $\,\pi^\prime
=\{N^\prime_1, \, \dots \,N^\prime_r\}$, $\sum_1^r N_j= \sum_1^r
N^\prime_j = N$ are equivalent if and only if the associated with
them functions $\psi_\pi :\{1, \, 2, \,\dots\, , \, N\} \to {\{0,
\,1,\, \dots\, r\}}$, and respectively $\psi_{\pi^\prime} :\{1, \,
2, \,\dots\, , \, N\} \to {\{0, \,1,\, \dots\, r\}}$ are equal.
Consequently  their  (common) orbit  is of length
$$  \big\vert \mathfrak{S}(r)/  \mathfrak{S}_\pi\big\vert = \frac{r !}{ \psi_\pi(1) ! \, \psi_\pi(2) !
    \cdots\,
\psi_\pi(N)!} $$
\end{teo}
\zz{Remind that the equivalence in $\Pi(r)$ means that there exists
a permutation $\sigma \in  \mathfrak{S}(r)$ acting on the indices of
partitions such that $ \pi^\prime = \sigma (\pi)$ (cf. Def.
\ref{definition equivalent partitions}).}

\subsection{The normalizer and Weyl group of a partial flag.}

\begin{defi}\label{normalizer of a flag}\rm
Let $\{V_j\}_1^r$ be an orthogonal partial flag. By the normalizer
of the flag in $O(N)$, respectively $SO(N)$, denoted
$\mathcal{N}(\{V_j\}_1^r)$, we mean the set of all elements $g$ of
$O(N)$, respectively $SO(N)$, which map $\{V_j\}_1^r$ into itself.
\end{defi}

Note that $\mathcal{N}(\{V_j\}_1^r)$ is a subgroup of $O(N)$,
respectively $SO(N)$,  containing the   orthogonal Borel subgroup,
respectively   connected orthogonal Borel subgroup of $\{V_j\}_1^r$.
Moreover, it is the normalizer of the latter in $O(N)$, respectively
$SO(N)$.

In this subsection we describe  the normalizer in $SO(N)$ of an
orthogonal partial flag $\{V_j\}_1^r$,
  or equivalently of an orthogonal   Borel
subgroup $H\equiv SO(N_1)\times \,\cdots \, SO(N_r)$, describe its
normalizer in $SO(N)$. This let us describe the set of all
orthogonal partial flags, or equivalently all orthogonal   Borel
subgroups $H$, for which the Weyl group
$\mathcal{W}(H)=\mathcal{N}(H)/H$ is nontrivial.

As a consequence,  we  characterize all    orthogonal Borel
subgroups which have a nontrivial Weyl group.

Let  $\{V_j\}_1^r$ be an orthogonal partial flag in $\br^N$ with
$\dim V_j =N_j$ giving a partition of $N$. Choosing an orthogonal
basis in each $V_j$ we get a coordinate system in $\br^N$ in which
$$\{V_j\}_1^r = \{\br^{N_1}, \br^{N_2}, \dots, \br^{N_r}\}\,.$$

\begin{prop}\label{Weyl group of flag}
Let  $\{V_j\}_1^r=\{\br^{N_1}, \br^{N_2}, \dots, \br^{N_r}\}$ be an
orthogonal partial flag in $\br^N$ with $\dim V_j =N_j$ defining a
partition $\pi$ of $N= \sum_{j=1}^r\, N_j$. Let $\pi=(N_1, \, \dots
\,, N_r)$ be a such that there is  $q$   different values $n_1,\,
\dots\,, n_q$ in the sequence $N_1, \, N_2,\, \dots\,, N_r$. Then
the normalizer $\mathcal{N}(\{V_j\}_1^r)$  in $O(N)$ contains the
orthogonal Borel subgroup $O(N_1)\times O(N_2)\times\, \cdots\,
\times O(N_r)$  and all permutations $\sigma_\pi \in
\mathfrak{S}_\pi =  \mathfrak{S}(\phi_\pi(n_1)) \times
\mathfrak{S}(\phi_\pi(n_2)) \times
    \cdots\, \times
\mathfrak{S}(\phi_\pi(n_q))\,.$ Moreover, the normalizer is
generated by these maps, and consequently the Weyl group of
$\{V_j\}_1^r$ is a product of $q$ permutation groups $
\mathcal{W}(\{V_j\}_1^r)= \mathfrak{S}_\pi =
\mathfrak{S}(\phi_\pi(n_1)) \times \mathfrak{S}(\phi_\pi(n_2))
\times
    \cdots\, \times
\mathfrak{S}(\phi_\pi(n_q))\,.$ Consequently, the Weyl  group
$\mathcal{W}(\{V_j\}_1^r)$ is of order $\psi_\pi(j_1) ! \,
\psi_\pi(j_2)! \cdots\, \psi_\pi(j_q)! = \psi_\pi(1) ! \,
\psi_\pi(2)! \cdots\, \psi_\pi(N)!\,,$ and it contains an element of
order $2$ provided it is not the trivial group.

\end{prop}
\noindent{\sc Proof.} It is clear that if an element  $g\in O(N)$
maps $\{V_j\}$ into itself then  for every $1\leq j \leq $ we have
$gV_j = V_{\sigma(j)}$ for some permutation $\sigma: \{1, \,
\dots\,, r\}$. Identifying every $V_j$
 with $\br^{N_j}$ by a choice of an orthonormal basis
$\{e_i^j\}$,  $1\leq i\leq N_j$, for any permutation $\sigma \in
\mathfrak{S}_r$ we can define a permutation map $g_\sigma \in O(N)$
by $$ g_\sigma(e_i^j) = e_i^{\sigma(j)}$$ provided only $
N_{\sigma(j)}=N_j$.

Composing $g$ with $g_{\sigma}^{-1}$ we get an orthogonal linear map
$ g^\prime  = g_{\sigma}^{-1}g: \br^N \to \br^N$ such that
$g^\prime(V_j) = V_j$ for every $1\leq j \leq  r$. Consequently
$g^\prime \in O(V_1)\times O(V_2)\times \, \cdots \, \times O(V_r)
\equiv O(N_1)\times O(N_2) \times\,\cdots\,\times O(N_r)$. This
shows that the normalizer $\mathcal{N}(\{V_j\})$ is generated by the
permutations $\sigma_\pi \in \mathfrak{S}_\pi$ and elements of its
  orthogonal Borel subgroup. \hbx

\begin{coro}\label{Weyl group of subgroup}
Let  $H \equiv O(N_1)\times O(N_2)\times\, \cdots\, \times O(N_r)$
be a     orthogonal Borel subgroup of $O(N)$. Then the Weyl group of
$H$ is equal to
 $   \mathcal{W}(H)=\mathfrak{S}_\pi =  \mathfrak{S}(\phi_\pi(n_1)) \times
 \mathfrak{S}(\phi_\pi(n_2)) \times
    \cdots\, \times
\mathfrak{S}(\phi_\pi(n_q)) $ $ =     \mathfrak{S}(\phi_\pi(1))
\times \mathfrak{S}(\phi_\pi(2))\times
    \cdots\, \times
\mathfrak{S}(\phi_\pi(r))  $ with a convention that
$\mathfrak{S}(\emptyset)= e$.
\end{coro}

\begin{rema}\label{combinatorics even}
Remind that with respect to our condition on the connected
orthogonal Borel
 subgroups  (cf. Definition \ref{Borel subgroup} and Remark \ref{dimension of conneceted Borel subroups}) in study of them we
 had to put the assumption $N_j \geq 2$. This  complicates a combinatorics  of possible partitions.
\end{rema}
\begin{rema}\label{Weyl of torus}\

It is worth of pointing out that if we take the partial orthogonal
flag $\{V_j\}$ in $\br^N$, $N=2K$, which is equivalent (i.e. is in
the same $O(N)$-orbit) to $\br^2 \oplus \br^2\, \oplus\, \cdots \,
\oplus \br^2$ then its   orthogonal Borel subgroup is conjugated to
$SO(2)\times SO(2)\times\, \cdots\, \times SO(2)= \bt^K$. In other
words it is the maximal  torus $\bt^K$ of $SO(2K)$.

On the other hand, by our considerations this conjugacy class of
Borel orthogonal subgroups
 corresponds to the partition $ 2+2+\,\cdots \, 2$ of $N=2K$,
Finally, from Corollary \ref{Weyl group of subgroup} it follows that
   $ {\text{for the  torus}} \;\; \bt^K \subset SO(2K)\;\;{\text{we have}}\;\;  \mathcal{W}(\bt^K)=
   \mathfrak{S}(K)\,,$
   which is a classical fact.
\end{rema}

\subsection{Partitions of $N$}

In the combinatorial number theory the number of ways of writing the
integer  as a sum of positive integers, where the order of addends
is not considered significant is denoted by $P(N)$, and    is
sometimes called {\it the number of unrestricted partitions}.

Also,   the number of ways of writing the integer  as a sum of
positive integers without regard to order with the constraint that
all integers in a given partition {\underline{are distinct}} is
denoted by $Q(N)$.

\begin{defi}[Definition of $P(N)$ and $Q(N)$]\label{definition of
P(N )}\rm  In our terms the function $P(N)$ is equal to the sum over
$1\leq r\leq N$ of numbers of equivalence classes of partitions $\pi
\in \Pi(r)$.

Respectively, the function $Q(N)$ corresponds to  the sum  over
$1\leq r\leq N$ of numbers of all equivalence classes of partitions
$\pi \in \Pi(r)$ such that all  $N_{j_1}, N_{j_2}, \, \dots\,,
N_{j_r}$ are different.

\end{defi}

More information about the functions $P(N)$ and $Q(N)$ we present in
Section \ref{Appendix}. At now let us only mention
  the asymptotic behavior of them (Hardy-Littlewood 1921):

 \begin{equation}\label{asymptotic}
 \begin{matrix}
 P(N)\,\sim\, \frac{1}{4 N\sqrt{3}} \, e^{\pi \sqrt{2N/3}}\,, \;\;\;\;
 &\;\;\;\;
 Q(N)\, \sim\,\frac{e^{\pi \sqrt{N/3}}}{4\, N^{3/}\,3^{1/4}}
 \end{matrix}
 \end{equation}

Let us remind that due to applications we are more interested in a
description of the number of conjugacy  classes of the
 orthogonal (connected) Borel subgroups of $O(N)$ ($SO(N)$)  with a
non-trivial Weyl group.

We denote by
\begin{equation}\label{definition of R(N)}  R(N):= P(N)-Q(N)
\end{equation}
 the number of
partitions of $N$ where the order of addends is not considered
significant, but which contain at least two equal summands.

As a consequence of Proposition \ref{Weyl group of flag} and
Corollary \ref{Weyl group of subgroup} we have the following.

\begin{teo}\label{number of classes}
Consider  the set $\mathfrak{F} ={\underset{1}{\overset{N}\cup}}
\mathfrak{F}_r$ of partial orthogonal flags in $\br^N$ of all
lengths $1 \leq r \leq N$ with the action of  $O(N)$, or $SO(N)$.

The number of classes of equivalence, i.e. the orbits of the action,
is  equal to $P(N)$. Consequently the number of  conjugacy classes
of all   orthogonal Borel subgroups of $O(N)$, or respectively
$SO(N)$ is equal to $P(N)$. Furthermore, the number of all
equivalence classes of partial orthogonal flags, respectively
orthogonal Borel subgroups, with the trivial Weyl group is equal to
$Q(N)$. Consequently the number of all equivalence classes of
partial orthogonal flags, respectively
  orthogonal Borel subgroups, with a non-trivial Weyl group is
equal to $\,R(N)= P(N)-Q(N)\,$.
\end{teo}

\begin{coro}\label{asymptotic of nontrivial Weyl}
The rate of growth of number $R(N)$ of the equivalence classes of
orthogonal partial flags of all lengths in $\br^N$ with a
non-trivial Weyl group is exponential:
$$
\begin{matrix}
 R(N) \sim P(N) - Q(N) \sim \frac{1}{4 N\sqrt{3}} \, e^{\pi
\sqrt{2N/3}}\,-\, \frac{e^{\pi \sqrt{N/3}}}{4\, N^{3/4}\,3^{1/4}}=
\,\frac{e^{\pi \sqrt{N/3}}}{4N^{3/4}} \, \Big(\frac{e^{\pi
(\sqrt{2N/3} - \sqrt{N/3})}}{ N^{1/4} \sqrt{3}} -
\frac{1}{3^{1/4}}\Big) \end{matrix}$$

\end{coro}

Due the analytical assumption  our task is to study equivalence
classes of partial orthogonal flags, respectively   orthogonal Borel
subgroups, such that each subspace is of dimension $\geq 2$.
 To do this we need a new notation.
\begin{defi}\label{partition without 1}\rm
For a given $N \in \bn$ let $P(N;1)$ denote the number of ways of
writing the integer $N$ as a sum of positive integers, where the
order of addends is not considered significant and each of them is
greater then $1$.

For a given $N \in \bn$ let $Q(N; 1)$ denote the number of ways of
writing the integer $N$  as a sum of positive integers without
regard to order with the constraint that all integers in a given
partition are distinct and greater then $1$.

Finally we define $R(N; 1):= P(N; 1)- Q(N; 1)$.
\end{defi}

We have the following correspondent of Theorem \ref{number of
classes} describing the number of classes  the   connected Borel
subgroups.

\begin{teo}\label{number of connected classes}
Consider  the set $\mathfrak{F} ={\underset{1}{\overset{N}\cup}}
\mathfrak{F}_r$ of partial orthogonal flags in $\br^N$ of all
lengths $1 \leq r \leq N$ with the action of $SO(N)$ and such that
$N_j=\dim V_j\geq 2$ for every summand of such a flag.

The number of classes of equivalence, i.e. the orbits of the action,
is  equal to $P(N; 1)$. Consequently the number of  conjugacy
classes of all   connected orthogonal Borel subgroups of  is equal
to $P(N;1)$. Furthermore, the number of all $O(N)$  equivalence
classes of partial orthogonal flags, respectively   connected
orthogonal Borel subgroups, with the trivial Weyl group is equal to
$Q(N; 1)$. Consequently the number of all equivalence classes of
partial orthogonal flags, respectively   connected orthogonal Borel
subgroups, with a non-trivial Weyl group is equal to $\,R(N; 1)\,$.
\end{teo}
\noindent{\sc Proof.} This  is a direct consequence of Proposition
\ref{Weyl group of flag}. Since in any such partition $N=N_1+N_2+\,
\cdots\, + N_r$ we have $N_j\geq 2$, we have take only such
partitions for which each summand is $>1$. The same with partitions
consisting of all distinct summands that correspond to flag, or
  connected orthogonal Borel subgroups, with the trivial Weyl
group. \hbx

Our next task is to give more  effective formulas for the the
functions $ \,P(N; 1)\,,Q(N; 1)$ and  $\,R(N; 1)$.

At this point we are {\bf not}   able to describe completely the
asymptotic behavior of the sequence $R(N;1)$.
 Regardless, we show that  there is  a family of  partitions of $N$ with each of them  with  nontrivial Weyl group and such   the rate of
growth of their amount   is exponential.

We define these partitions by a formula which has four different
forms depending on the class of $N$ modulo $4$.

Suppose first that  $N\cong 0 \mod 4$,  is  i.e. it is of the form
$N=4M$, $M\geq 1$. Let $ M_1 \leq M_2 \leq \, \dots\, \leq  M_r$ be
a partition of $M$ od length $r$.

We define  a partition of $N=4M$ of length $ r$  by the formula
\begin{equation}\label{first double partition}
\pi_r(N)= \{2\,M_1, \, 2\,M_1, \, 2M_2,\, 2M_2,  \,\dots \,, 2M_r,\,
2\,M_r\}
\end{equation}
For every $1\leq j\leq r$ we have $N_{2j}=N_{2j+1} = 2M_j \geq 2$
and the Weyl group  $\mathcal{W}(\pi_r(N))$ is nontrivial. Indeed
the Weyl group of contains the Weyl group of the partition $ M_1
\leq M_2 \leq \, \dots\, \leq M_r$ of $M$ and transpositions of the
equal summands $2 M_i \mapsto 2 M_i$.

Suppose now that  $N \cong 2 \mod 4$, $N\geq 6$, i.e. $N$  is   of
the form $4M +2$, $M\geq 1$. Let $ M_1 \leq M_2 \leq \, \dots\,
,\leq M_r$ be a partition of $M$ of length $r$. We define a
partition of $N=2M+2$ of length $r$ by the formula

\begin{equation}\label{second double partition}
\pi_r(N)=  \{2, 2\,M_1, \, 2\,M_1,\, 2\,M_2,\,2\,M_2, \,\dots\,, 2\,
M_r  ,\, 2\,M_r \}
\end{equation}

For every $1\leq j\leq r$ we have $N_j \geq 2$ and the Weyl group
$\mathcal{W}(\pi_r(N))$ is nontrivial. Indeed   the Weyl group of it
contains the Weyl group of the partition $ \{ M_1 \leq M_2 \leq \,
\dots\, \leq M_r\} $ of $M$ and transpositions of the equal summands
$2M_i \mapsto 2 M_i$ for $i\geq 2$.

Thirdly, suppose that  now  $N\cong 3 \mod 4 $,   $N\geq 7$, i.e.
$N$ is of the form $N=4M+3$, $M\geq 1$.  Let $ M_1 \leq M_2 \leq \,
\dots\, ,\leq M_r$ be a partition of $M$ of length $r$. We define a
partition of $N=4M+3$ of length $r$ by the formula

\begin{equation}\label{third double partition}
\pi_r(N)= \begin{cases} \{2\,M_1, \, 2\,M_1,\,2\,M_2,\, 2\,M_2,\,
\dots \,,2 M_i,\, 2\,M_i,\, 3,  \,2 M_{i+1}, \,2\,M_{i+1},\, \dots
\,, 2\,M_r \} \;{\text{if}}\;\; M_1, M_2,\dots\,, M_i \leq 1 \,, \cr
\{3, \, 2\, M_1,\, 2\, M_2,\,\dots\,, 2\,M_r\} \;{\text{if}}\;\;
2\leq M_1\leq M_2\leq\, \cdots \, \leq M_r\,,
\end{cases}
\end{equation}
For every $1\leq j\leq r$ we have $N_j \geq 2$ and the Weyl group
$\mathcal{W}(\pi_r(N))$ is nontrivial by the same argument as above.

Finally, suppose that $N \cong 3 \mod 4 $, $N\geq 9$, i.e. $N=
4M^\prime +1$, $M^\prime \geq 2$, or equivalently $N= 4M +5$, where
$M=M^\prime -1$. Let $ M_1 \leq M_2 \leq \, \dots\, ,\leq M_r$ be a
partition of $M$ of length $r$. We define a partition of $N=4M+5$ of
length $r$ by the formula
\begin{equation}\label{forth double partition}
\pi_r(N)= \begin{cases} \{2\,M_1, \, 2\,M_1,\,2\,M_2,\, 2\,M_2,\,
\dots \,,2 M_i,\, 2\,M_i,\, 5,  \,2 M_{i+1}, \,2\,M_{i+1},\, \dots
\,, 2\,M_r \} \;{\text{if}}\;\; M_1, M_2,\dots\,, M_i \leq 2 \,, \cr
\{5, \, 2\, M_1,\, 2\, M_2,\,\dots\,, 2\,M_r\} \;{\text{if}}\;\;
3\leq M_1\leq M_2\leq\, \cdots \, \leq M_r\,,
\end{cases}
\end{equation}

For every $1\leq j\leq r$ we have $N_j \geq 2$ and the Weyl group
$\mathcal{W}(\pi_r(N))$ is nontrivial by the same argument as above.

 As a consequence we get the following.
\begin{prop}\label{special partitions}
For every   $N\geq 4$ of the form $N=4M$,   $N=4M +2$, $N=4M+3$, or
finally $N=4M+5$ with $M\geq 1$ there exists at least $s_N=
P(\frac{N}{4})$, $s_N= P(\frac{N-2}{4})$, $ s_N= P(\frac{N-3}{4})$,
or respectively $s_N=P(\frac{N-5}{4})$
 equivalence classes of partitions of $N$
such that every summand $N_j\geq 2 $ and for every such partition
$\pi(N)$ the Weyl group is nontrivial.
\end{prop}

\section{Applications to analytical problems}
For  analytical applications  we have to show that for any  two not
equivalent   orthogonal  Borel  subgroups $H,\,K \subset O(N)$ the
subgroup $\langle H, K \rangle $ generated by them acts transitively
on a sphere $S(V)$, where $V \subset \br^N$ is a subspace invariant
with respect to  $H$ and $K$ simultaneously. The resultating answer
says that it does not depend on a choice of a representatives of the
conjugacy classes of $H$, and $K$ in $O(N)$.

More precisely, by $G=\langle H, K\rangle$ we denote a subgroup of
$O(N)$ {\it generated topologically} by $H$ and $K$, i.e. the
closure of the set of products $\{h^{a_1} k^{b_1}\, \dots\, h^{a_m}
k^{b_m}\}$, where $h\in H, \,k\in K$, $a_i\in \mathbb{Z}, \; b_i\in
\mathbb{Z}$.

\begin{teo}\label{main}
Let $H,\, K \subset O(N)$ be two  Borel subgroups

 If $\langle H,K
\rangle $ acts transitively on $S(\br^N)$ then for every $g\in O(N)$
the group $\langle H, gKg^{-1} \rangle $, and dually $ \langle
gHg^{-1}, K \rangle $ acts
 transitively on $S(\br^N)$.

 Moreover, suppose that there exists $N^\prime < N$ such that
 $H(\br^{N^\prime})
 \subset \br^{N^\prime}$, $K(\br^{N^\prime})
 \subset \br^{N^\prime}$ then  $ \langle H, K \rangle  (\br^{N^\prime}) \subset
 \br^{N^\prime}$ and  $  \langle H, K \rangle$ acts transitively on the sphere $S(\br^{N^\prime})$.

 Then statement holds  for every $g^{\prime} \in
 O(N^\prime) \subset O(N)$, i.e.
 the group $ \langle H, g^\prime K {g^{\prime}}^{-1} \rangle$, and dually $ \langle g^\prime H{g^\prime}^{-1}, K \rangle $
acts
 transitively on $S(\br^{N^\prime})$.

\end{teo}

\begin{coro}
$\langle H,K\rangle $ acts transitively on $S(\br^N)$ if and only if
for every
 $g, \,\tilde{g} \in O(N)$ the subgroup $\langle gHg^{-1},
 \tilde{g}K\tilde{g}^{-1}\rangle $ acts transitively on $S(\br^n)$.
\end{coro}

\noindent \textsc{Proof of Theorem \ref{main}.} We begin with
proving the first part of Proposition \ref{main}, i.e. that the
transitive action of  $\langle H,K\rangle $  on $S(\br^N) $ implies
the transitive action of $\langle H, gKg^{-1}\rangle $ on $S(\br^N)
$.

First of all note due to Remark \ref{action of component of
identity} for $N\geq 2$ $\langle H, K\rangle$ acts transitively on
$S(\br^N)$  if and only if  $\langle H, K\rangle_0$ does.
Consequently, $\langle H, K\rangle$ acts transitively on $S(\br^N)$
if and only if  $\langle H_0, K_0\rangle$ does. From it follows that
we can assume that $H$ and $K$ are connected.

Next, we assume that $g\in SO(N)\subset O(N)$ leaving the remaining
case to the end of proof.

 Now we show
that the statement holds if $g$ is in a small neighborhood of $e\in
SO(N)$. Denote $\langle H, K\rangle $ by $G\subset SO(N)$ and
$S(\br^N)$ by $\mathcal{M}$. The transitivity of action of $ G
\subset SO(N)$   means that  the map $ \phi_m: SO(N) \times \{m\}
\to \mathcal{M}$ restricted to $G\times \{ m\} $ is a surjection. It
is enough to show that locally at $ (e, m)$ it is a surjection. But
this  is equivalent to the fact that the linear map $D\phi(e,m) :
TG_e \to TM_m$ is a surjection.

The latter is equivalent to the fact that ${\rm rank} \,
D\phi(e,m)_{| T_e G}=  \dim T_m S^{N-1} = N-1$. Now observe that the
subspaces $T_e G$ and $T_eG^\prime$,  where $G^\prime = \langle H,
gKg^{-1}\rangle $, are close each to the other if $g$ is close to
$e$. This means that they are close as elements of the Grassmanian
$G(l, d)$,  where $l=\dim SO(N)$ and $d=\dim G=\dim G^\prime$.
Equivalently the distance:
$$ \max \rho(x,y): \{x\in T_e G\cap S(\br^N), \, y\in  T_eG^\prime
\cap S(\br^N)\} $$ is small. If we we change   $H$ to   $K$, and
conversely, then the proof is the same.

 Now the statement follows from the  continuity of
$D\phi(e,m)$. Indeed, for a basis $\{x_1, \, \dots\, x_l\}$, $ x_i
\in T_eG^\prime \cap S(\br^N)$  be a basis of $T_e G$ such that
${\rm rank}  \{D\phi(e,m)(x_i)\} = N-1$. By the continuity of
$D\phi(e,m)$ if $\{x_1^\prime, \, \dots\, x_l^\prime\}$ is a basis
of $T_e G^\prime$ such that $ \max \rho(x_i,x_i^\prime)$ is small,
then
$${\rm rank}  \{D\phi(e,m)(x_i^\prime)\} ={\rm rank}
\{D\phi(e,m)(x_i)\} = N-1$$, which shows that
 $ D\phi(e,m): T_e G^\prime \to  T S^N$ is a surjection if $g$ is close
 to $e$.

Now let $\tilde{g} = g_0 g$ be an arbitrary element of $O(N)$ with
$g\in SO(N)$ and $g_0 \notin SO(N)$.  We can assume that $g_0$ is an
involution $g_0^2 = {\rm id}$ with $\det (g_0)= -1$. It is clear
that $\langle H, g_0 (g Kg^{-1}) g_0^{-1}\rangle$ acts transitively
on the sphere $S(V)$ if and only if $\langle H,g Kg^{-1}) \rangle$
does so, which reduces this case to the previous.

  \vskip 0.1cm

 Now let us take $H= gHg^{-1}$ and $K=
g(gKg^{-1})g^{-1}= g^2 K (g^2)^{-1}$.

By the above $\langle H,K\rangle $, $\langle H, gKg^{-1}\rangle $,
$\langle gHg^{-1}, K\rangle $, and $\langle gHg^{-1},
gKg^{-1}\rangle $ act transitively on $S(\br^N)$. The latter is
obvious, because $\langle gHg^{-1}, gKg^{-1}\rangle  = g\langle
H,K\rangle  g^{-1} $. This implies that $\langle g Hg^{-1},g^2 K
(g^2)^{-1}\rangle $ acts transitively on $S(\br^N)$. Now applying
the first part of this  proof to $H$, with $g=g^{-1}$ being small,
we see that $\langle H, g^2Kg(^2)^{-1}\rangle $ acts transitively on
$S(\br^N)$.

Continuing this argument, we see that for every $n\in \bn $ $\langle
H, g^n K (g^n)^{-1}\rangle $, and  respectively $\langle g^n H
(g^n)^{-1}, K \rangle $ act transitively on $S(\br^N)$ if $g\in U$
belongs to some small symmetric (i.e. $U= U^{-1}$) neighborhood of
$e$. But it is known that any such neighborhood of $e$ generates
every connected compact Lie group, which completes the proof of
first part of statement.

A proof of the second part  with a  $ V^\prime\varsubsetneq V$
invariant for $H$ and $K$  is analogous.  \hbx

\vskip 0.3cm

The above theorem let state that  the  conjugacy classes of the
  Borel subgroups are determined by partitions of
the number $N= \sum_j^r N_j$ and isomorphism classes of Borel groups
of dimensions $N_j$. In particular, the  conjugacy classes of the
   orthogonal  Borel subgroups are determined by the equivalency classes of partitions of
the number $N= \sum_j^r N_j$.
\begin{lema}\label{determination by partition}
Let $ H =  H_1\times H_2 \times \, \cdots\, \times H_r$
 where $H_j = \pi_j(H)$, $H_j \subset O(N_j)$, and ${\overset{r}{\underset{1}\sum}} N_j=N$ be an
 orthogonal Borel subgroup  corresponding to a partial flag
 $\{0\}= V^0 \varsubsetneq V^1 \varsubsetneq V^2\, \dots
\varsubsetneq \, V^r=V\;\;{\text{in}} \;\; V=\br^N $  with its
orthogonal decomposition $ V={\underset{j=1}{\overset{r} \oplus}}
V_j, \; \; {\text{where}}\;\; V_j = V^{j-1}_\perp \supset V^j  $
spanned  by an orthonormal basis $\{e^i_j\}$, $1\leq j \leq r$,
$1\leq i\leq N_j$ of $\br^N$. Let next $ K =   H_1^\prime\times
H_2^\prime \times \, \cdots\, \times H_s^\prime$ be another
orthogonal subgroup,
 where $H_j^\prime = \pi_j(K)$, $H_j^\prime \subset O(N_j^\prime)$, and ${\overset{s}{\underset{1}\sum}} N_j^\prime=N$ be an
 orthogonal Borel subgroup  corresponding to a partial flag
 $\{0\}= {V^\prime}^0 \varsubsetneq {V^\prime}^1 \varsubsetneq {V^\prime}^2\, \dots
\varsubsetneq \, {V^\prime}^s=V\;\;{\text{in}} \;\; V=\br^N $  with
its orthogonal decomposition $ V={\underset{j=1}{\overset{s}
\oplus}} {V^\prime}_j, \; \; {\text{where}}\;\; V_j^\prime =
{V^\prime}^{j-1}_\perp \supset {V^\prime}^j  $ of $\br^N$  spanned
by an orthonormal basis $\{{e^\prime}^i_j\}$, $1\leq j \leq s$,
$1\leq i\leq N^\prime_j$.

Then there exists a representative $\tilde{K} = g Kg^{-1}$ of the
conjugacy class of $K$ such that the corresponding partial flag
$\{\tilde{V}^j \}$, $\tilde{V}^j =
{\underset{j=1}{\overset{s}\oplus}} \tilde{V}_j$, $1\leq j \leq s$
is spanned by the  basis $\{e^i_j\}$ which spans  the partial flag
$\{V^j\}$ corresponding to the subgroup $H$.
 Moreover, if
$\{e_j^i\}$, and $\{{e_j^\prime}^i\}$ are of the same orientation,
in particular  if $H$ and $K$ are connected orthogonal Borel
subgroups, then $g \in SO(N)$. Furthermore $\tilde{N}_j =
N_j^\prime$ for every $ 1\leq j \leq s$.
\end{lema}
\noindent \textsc{Proof.} By the linear algebra we know that
ordering ${e^\prime}^i_j \mapsto e_j^i $, $1\leq i\leq N$ extends to
an orthogonal map $ g\in O(N)$. If $\{e^i_j\}$  and
$\{{e^\prime}^i_j\}$ are in the same orientation class then $g \in
SO(N)$. Then $ \tilde{V_j}=g (V_j^\prime)$, and respectively $
\tilde{V^j}= g (V_j^\prime)$ are the required partial flags of the
statement.

Since $g$ is an isomorphism $ \tilde{N}_j= \dim \tilde{V}_j = \dim
V^\prime_j= N^\prime_j $  which completes the proof.\hbx

As we already noted,  in the basis $\{e^i_j\}$ the orthogonal Borel
subgroup $H$, and correspondingly the subgroup $\tilde{K}$, defines
a partition  $N_j, 1\leq j \leq r$ of $N$, respectively another
partition $\tilde{N}_j= N^\prime_j$, $1\leq j \leq s$ of $N$.
Consequently, in view of Theorem \ref{main} and Lemma
\ref{determination by partition}  to  study of the problem of action
of group $\langle H, K\rangle $ on the sphere we can take the
representatives of $H$ and $K$ in the same basis.

\zz{The Lemma  \ref{determination by partition} shows that   the
maximal orthogonal Borel subgroups associated with a given fixed
orthonormal basis of $\br^N$ give  all equivalence classes of  pairs
$H, K$ of possible maximal orthogonal groups Borel subgroups such
that $\langleH, K \rangle$ acts transitively on}

\zz{Note that any orthogonal map $g: \br^N\to \br^N $ given by a
permutation of  the basis $\{e_j^i\}$  in such a way that the
induced orthogonal map $g$ permutes subspaces $\{V_j\}$ $ 1\leq
j\leq r$, and respectively that the subgroup $gHg^{-1}$ respectively
a subspaces $\{\tilde{V}_j\}$ $ 1\leq j\leq s$, induces partitions
of $N$ which re permutations of original, thus equivalent to the
original.}

This let us adapt and apply the geometrical result of \cite[Lemma
4.1]{KrMa}.

Let $N^H_1 \leq N^H_2\leq N^H_3\leq \, \cdots\, \leq N^H_r$, $N_1^K
\leq N^K_2\leq \, \cdots\, N^K_s$, $\sum_1^r N_j^H = \sum_1^s N^K_i
=N$ be two different  not  decreasing partitions  of $N$. Let $1\leq
a  \leq r,  \, s$,  $\,N_a\leq N$,  be the smallest numbers such
that for all   $1\leq j \leq a $, $1\leq i\leq a $ we have   $N^H_j
= N^K_i$, and  $$ N_{a} = \sum^{a}_1 N^H_i =\sum_1^{a}N_i^K$$ with
the convention that a sum over empty set is zero, i.e. in the case
when there is not such a  $a <r, \, a < s$. In other words, up to
the index $a$ these partitions are equal and define a partition of
$N_{a}\leq N$. If $a=r=s$ then these partitions of $N$ are equal.

Let next $ N_{b} \leq  N $, $ N_b >N_{a}$ be a smallest number such
that there exist $
 r\geq j_{b} >  a$ and $ s\geq  i_{b} >i_a$ such that
$$ N_a + \sum_{j_a+1}^{j_b} N^H_j =  N_a+ \sum_{i_a+1}^{i_b} N^K_i  = N_b$$
Of course, it can happen that $j_b=r$. Then $i_b=s$, and
consequently $N_b=N$.

Put $N_{{b}-{a}}:= N_{b}- N_{a}$. By the above $N^H_{j_a+1},
N^H_{j_{a_1}+2},\, \dots \, N^H_{j_{b}}$ and $N^K_{i_{a}+1},
N^H_{i_a+2},\, \dots \, N^H_{i_{b}}$  are partitions of $N_{b-a}$
and they do not contain  a common  sub-partition of any $N^\prime <
N_{b-a}$. As a consequence of \cite[Lemma 4.1]{KrMa} we have the
following proposition

\begin{prop}\label{adaptation of KrMa}
Let $N^H_1 \leq N^H_2\leq N^H_3\leq \, \cdots\, \leq N^H_r$, $N_1^K
\leq N^K_2\leq \, \cdots\, N^K_s$, $\sum_1^r N_j^H = \sum_1^s N^K_i
=N$ be two different  not  decreasing partitions  of $N$. Let $a,
\,b$ , correspondingly $N_a < N_b$ , and $N_{b-a} \leq N$ be as
defined above. Let next $H= O(N_1^H) \times O(N^H_2) \times \,
\cdots\times O(N^H_s)$, correspondingly $K= O(N^K_1)\times O(N^K_2)
\times \,\cdots\, \times O(N^K_s)$, be the maximal orthogonal Borel
subgroups associated with these partitions.

Then the subspace $V^\prime = \{0\}\times \{0\}\times \cdots  \times
\br^{N^H_{j_a+1}} \times \br^{N^H_{j_a+2}}\times \cdots \times
\br^{N^H_{j_b}}\times \{0\}\times \cdots \{0\}$  of dimension
$N_{b-a}$ is preserved by $H$ and $K$ and the group $\langle H, K
\rangle$ acts transitively on the sphere $S(V^\prime)$.
\end{prop}

Our main analytical observation is the fact that to apply the scheme
of \cite{KrMa} that is based on \cite{BaWi} we do not need use pairs
$H,\, K$  of subgroups, with nontrivial the Weyl group, which
generate an entire Borel group $ G=\langle H, K\rangle \subset
O(N)$. It is enough if they generate a group which acts transitively
on $S(V^\prime)$  for a subspace $V^\prime \subset \br^N$ preserved
by both $H$, and $K$.

Let  $E$ be  a functional space on which $H, \, K \subset O(N)$ act
by the action induced by the action of $O(N)$ on the space of
variables in $\br^N$. Suppose the $\tau \in \mathcal{N}(H)\setminus
H$, and correspondingly,  $ \tau^\prime \in \mathcal{N}(K)\setminus
K$ are elements of order two in $\mathcal{N}(H)$, and
$\mathcal{N}(K)$ respectively. Suppose next that $\Upsilon
=\{\tau_\alpha\}$, and $\Upsilon^\prime=\{\tau^\prime_\beta\}$ are
finite sets of such elements in $\mathcal{N}(H)$, and
$\mathcal{N}(K)$ respectively.

Consider  $ \tilde{H}= \langle H, \,\Upsilon\rangle $, and
correspondingly $ \tilde{K}= \langle K, \,\Upsilon^\prime \rangle $,
 the extended subgroups of $\mathcal{N}(H)$ generated by $H$ and
$\Upsilon$, and $K$ and $\Upsilon^\prime$ respectively.

 Now define an action  of the group $\tilde{H}\subset \mathcal{N}(H)$ on
 $E$,
 and correspondingly an action of   $\tilde{K}\subset \mathcal{N}(K)$, as a representation) $\rho :
\tilde{H}  \to GL(E)$, respectively $\rho : \tilde{H}  \to GL(E)$.
 For an element of the form
$g= \tau h$, $\tau \in\Upsilon$, $h\in H$, or   respectively
$g^\prime=\tau^\prime k$, $\tau^\prime \in\Upsilon$,   $k\in K$, it
is given by the formula:
$$  \rho(h) f (x) : =   f(h  x)
{ \;\;\;}  \rho (\tau) f (x)  = -1 f( \tau x) \,,
$$
and next extended to a homomorphism $\tilde{H}\to GL(E)$. Obviously,
it extends to a homomorphism from $\tilde{H}$ to $GL(E)$ and does
not depend on the representative  of $g$ as $ \tau h$, since
$\tau\in \mathcal{N}(H)$ and $-{\rm Id}$ commutes with every linear
map of $E$. Similarly we define $\rho: \tilde{K}\to GL(E)$.

We will consider  the subspaces $E^{\tilde{H}}$, respectively
$E^{\tilde{K}}$,   of the fixed points sets of the above defined
action.

he following observation is fundamental
\begin{prop}\label{linear independence of fixed points}
Let $E$ be a function space with an $O(N)$-invariant  domain, $N\geq
4$. Let next $H$, $K$ be two subgroups of $O(N)$ such that:
\begin{itemize}
\item{the Weyl groups $\mathcal{W}(H)$, $\mathcal{W}(K)$ are nontrivial and contain elements of order two $\tau, \,\tau^\prime$ respectively.}
\item{ there exists $2\leq N^\prime \leq N$ such that the group
$G=\langle H, K\rangle $ acts transitively on $S(V^\prime)$ for
subspace $V^\prime \simeq \br^{N^\prime}$ preserved by $H$, and
$K$.}
\item{there exist  $\tau \in \Upsilon$ or  $\tau^\prime \in
\Upsilon^\prime$ which acts  trivially on  the orthogonal complement
$V^\prime_\perp$ and acts nontrivially on
 $V^\prime$.}
\end{itemize}
Then for the groups $\tilde{H}, \tilde{K}$ defined above we have
$$ E^{\tilde{H}} \cap E^{\tilde{K}} =\{0\}\,.$$

\end{prop}
\noindent{\sc Proof.} First note that if $H(V^\prime)\subset
V^\prime$ and $K(V^\prime)\subset  V^\prime$ then also
$H(V^\prime_\perp)\subset V^\prime_\perp$ and respectively
$K(V^\prime_\perp)\subset V^\prime_\perp$. Let $(x,y)$  be the
coordinates in $\br^N$ corresponding to the orthogonal decomposition
$\br^N = V^\prime \oplus  V^\prime_\perp$.

Assume that there exists $\tau \in \Upsilon$  such that $\tau$ acts
nontrivially on $V^\prime$ and trivially on $ V^\prime_\perp$.

Let $f\in E^{\tilde{H}} \cap E^{\tilde{K}} \subset E^H\cap E^K$.
Take $\bar{H}=\langle H, \tau \rangle \subset \tilde{H}$. So that we
have  $f\in E^{\tilde{H}} \cap E^{\tilde{K}} \subset E^{\bar{H}}
\cap E^{\tilde{K}}\subset E^H\cap E^K$.

Fix $y \in V^\prime_\perp$. Since $G=\langle H, K\rangle $ acts
transitively on $S(V^\prime)$, the function $f$ is radial in $x$,
i.e. $f(x, y)$ depends only on $\vert x \vert$ as $ g f(x,y) = f(gx,
gy)$ for every $ g \in \langle H, K\rangle$ and $\langle H,
K\rangle$ acts transitively on $S(V^\prime)$.

  Moreover
$f(\tau   (x,y)) = f(\tau \,x, y)$, because $\tau$  preserves
$V^\prime$, and acts trivially on $V^\prime_\perp$. This gives
$$ -f(x,y)  =   f(\tau (x,y)) = f (\tau  x, y) =  f(x,y)\,, $$ for every
$x \in  V^\prime$, since $\vert \tau  x \vert = \vert x \vert $.

This shows that for every fixed $y$ the section function $f_y(x)=
f(x,y)$ is identically equal to $0$. Consequently $f(x,y)$ is equal
identically zero which proves the statement. \hbx

As we announced we are not going into  analytical complexity of
particular variational problems that need a special study by fine
analytical tools. Instead  we assume what usually is  an output of
this study as our supposition (see \cite[Th. 3.2]{BaWi}).

\begin{assum}[{\bf  On the invariant functional}
$\Phi$]\label{Assumption on functional} {\;\;}

 Let $E$ be an infinite
Hilbert space with an orthogonal linear action $\rho(g): E \to E$ of
a compact Lie group $G$, $\Phi: E\to \br$ is a $C^1$ functional and
the following hold:
\begin{itemize}

\item[i)]{
$\Phi$ is $G$-invariant.}
\item[ii)] {$\dim E^G = \infty$.}
\item[iii)] {$\Phi_{|E^G}: E^G \to \br$  has infinitely many critical
points  $\{u_k\}$ such that $\Vert u_k\Vert \geq c_k$, $ c_k \to
\infty$.}

\end{itemize}
\end{assum}

\begin{rema}\label{explaining of dimension} In all the quoted papers \cite{BaWi},
\cite{BaLiWe3}, \cite{KrMa}, \cite{Kr} and others, e.g. also quoted
in \cite{KrRaVa}, the Assumption \ref{Assumption on functional} is
verified for studied there problems. It is based on a scheme  which
is called the fountain theorem (see \cite{Ba} for an excellent
exposition), and which began with the famous Ambrosetti Rabinowitz
theorem of \cite{AmRa}. The latter is applicable, because  in all
the cited above works the studied problems lead to a functional
$\Phi$ which is even, i.e $\Phi(-u) = \Phi(u)$.

Notify, since the functional space $E$ consists  of functions with
domain  equal to $\br^N$, the Palais-Smale condition is not
satisfied for a functional $\Phi: E \to  \br$ in general. But her
$\Phi$ is $O(N)$-invariant with the action of $O(N)$ on the
functional space $E$ induced by the action of $O(N)$ on  the domain
of functions equal to $\br^N$, and $G=\langle H, \tau\rangle$,  $H =
H_1\times H_2\,\cdots\, H_r$, $H_j \subset O(N_j)$, $\tau \in
\mathcal{W}(H)\setminus H$,  $\sum N_j= N$. Consequently,   the
Palais-Smale condition follow form the Lions theorem \cite{Lions}
which says that for the discussed functional space the embedding
$$ E^G \hookrightarrow L^s(\br^N)$$
is compact for  $s$ of an interval $[a,b]$ depending on the problem
(on $N$) provided $ N_j \geq 2$ for every $1\leq j \leq r$.

This additionally to our Remark \ref{action of component of
identity} justifies the assumption $ N_j \geq 2$ put on an
orthogonal Borel subgroup.
\end{rema}

In the terms used here, the  problem of multiplicity of infinite
series of non-radial sign-changing solutions was reduced to a
problem of finding a  set $\{H_i\}_{i=1}^{s_N}$, with some $s_N$
depending on $N$, of subgroups of $ O(N)$
such that: \\

- for every $i\neq j$ $H_i$ and $H_j$ are not conjugated in $O(N)$,

- for every $i\neq j$ the group $\langle H_i, H_j\rangle$ acts
transitively on $S(V^\prime))$, where  $V^\prime \simeq \br^N$ and
is invariant for $H$ and $K$,

- for each $1\leq i\leq s_N$ there exists an element  $e\neq \tau
\in \mathcal{W}(H_i)$  of order $2$, or   such that
 $\tau(V^\prime) \subset V^\prime$,  and $ \tau_{|V^\prime_\perp}= {\rm id}_{|V^\prime_\perp}$, or an element  $e\neq
 \tau^\prime
\in \mathcal{W}(H_i)$  of order $2$ of the same property. \\

With the notation of previous section and Theorem \ref{number of
classes} we have the following theorem.

\begin{teo}\label{final theorem}
Let $N \geq 4$ and $O(N)$ orthogonal group of $\br^N$ and $\Omega
\subset \br^N$ be equal to $\br^N$, $D_r^N$ or $S^{N-1}$. Suppose
that we have a variational  problem which weak  solutions correspond
to the critical points of a functional $\Phi: E \to  \br $ defined
on a functional space $E$ on  $\Omega$, where $E$ posses the natural
linear action induced by the action of $O(N)$ on $\br^N$.
 Assume also that $\Phi$ is $O(N)$-invariant and satisfies
 Assumption \ref{Assumption on functional} with respect to  every
 closed subgroup $G\subset O(N)$.

Then there exist $s_N$ geometrically distinct, with respect to the
symmetry, infinite series of solutions $ \mathcal{S}_i =\{u^i_k\}$,
$1\leq i \leq s_N$, $1\leq k < \infty $ of this problem, where
$$ s_N\,=\begin{cases} P(\frac{N}{4}) \;\;{\text{if}} \;\; N\geq 4 \;\;{\text{and}}\;\;  N=4M\,,\cr
P(\frac{N-5}{4})\;\;{\text{if}} \;\; N\geq 9 \;\;{\text{and}}\;\;
N=4M +5 = 4M^\prime +1 \,\cr
 P(\frac{N-2}{4})\;\;{\text{if}} \;\; N\geq 6
\;\;{\text{and}}\;\; N= 4M +2 \,\cr
P(\frac{N-3}{4}) \;\;{\text{if}}
\;\; N\geq 7 \;\;{\text{and}}\;\;N=4M +3\,.
\end{cases}
$$
with $P(N)$ defined in (\ref{definition of P(N )}) with an
exponential asymptotic behavior given in (\ref{asymptotic}).
Moreover,  in each series $\mathcal{S}_i$ the solutions $u_k^i$ are
geometrically distinct for different $k$.

\end{teo}
Note that Theorem \ref{final theorem} applies to the problems
studied in \cite{BaWi}, \cite{BaLiWe3}, \cite{KrMa}, and
\cite{KrRaVa}.. Another trivial observation: $N=5$ is not covered by
the cases listed above. Any other number $N\geq 4$ appears there.

\noindent {\sc{Proof.}} We will prove the theorem  showing that the
assumptions of Proposition \ref{linear independence of fixed points}
are satisfied for the groups determined by the partitions given in
the statement. To do it we first have to show that the assumption of
Proposition \ref{adaptation of KrMa} is satisfied.

 First of all let
us take the unique representative $M_1 \leq M_2\leq M_3 \leq
\,\cdots \, \leq M_r$, of each equivalence class of $P(M)$
partitions of  $M= \sum M_j$ \zz{with nontrivial Weyl group}. Then
the formula defined in (\ref{first double partition}),(\ref{second
double partition}), (\ref{third double partition}), or respectively
(\ref{forth double partition}) gives a partition of $N$ in each of
listed above cases  modulo $4$. Obviously  for every summand we have
$N_j\geq 2$.

Assume for simplicity that $N=4M$, i.e we have the first of four
cases. Let $\pi(M)$ and $\pi^\prime(M)$ be two different partitions
of $M$. Observe that if $0\leq M_a < M_b \leq M$ is the interval on
which $\pi(M)$ and $\pi^\prime(M)$ do not coincide, with the minimal
$M_a$ with this property, then an interval $0\leq N_a = 4M_a< N_b
\leq 4M_b\leq 4M$ has the corresponding property with respect to the
partitions of $N=4M$:

\begin{equation}\label{double partitions}
\begin{matrix}
 \{2M_1, 2M_1, 2 M_2, 2M_2,\, \dots\,,2M_{j_a}, 2M_{j_a},\, 2M_{j_a+1}, 2M_{j_a+1} ,
\, \dots\,,, 2M_r, 2M_r\} \cr {\; \; } \cr \{2M_1, 2M_1, 2 M_2,
2M_2,\, \dots\,,2M_{i_a}, 2M_{i_a},\, 2M^\prime_{i_a+1},
2M^\prime_{i_a+1}, \, \dots\,,, 2M^\prime_s,\,2M^\prime_s\}\,.
\end{matrix}
\end{equation}

Indeed $M_j=M^\prime_i$ for $i=j\leq j_a=i_a$, so $2M_j=2M^\prime_i$
for $i=j\leq j_a=i_a$, $2M_{j_a+1} \neq  2 M^\prime_{i_a+1}$ as
$M_{j_a+1} \neq  M^\prime_{i_a+1}$ thus $2j_a=2i_a$ is the maximal
index up to which partitions (\ref{double partitions}) of $N=4N$
coincide.

We have  $N_b \leq 4M_b$ but it could happen that  $N_b < 4M_b$ in
general. Nevertheless, the interval $N_a, N_b$  on which the
partitions do not coincide contains at least one pair $M_{j_a+1}, \,
M_{j_a+1}$, or dually $M^\prime_{i_a+1}, \, M_{i_a+1}$, because it
can happen that $ 2 M_{j_a+1} + 2M_{j_a+1} = 2 M^\prime_{i_a+1}$ or
dually $ 2 M^\prime_{i_a+1} + 2 M^\prime_{i_a+1}= 2 M_{j_a+1}$.

Now let us take  as $N_b$ the largest  number such that the
partitions $\pi(N)$ and $\pi^\prime(N)$ do not coincide in any
subinterval, and $\tau$, or $\tau^\prime$, the transposition of
first pair of coordinates corresponding to $2 M_{j_a+1},2
M_{j_a+1}$, ro respectively $2 M^\prime_{i_a+1}, 2 M^\prime_{i_a+1}$
depending which of the above cases happens. (Note that the both
cases  could happen). Let $\tilde{j}_b$, or dually $\tilde{i}_b$ be
the index corresponding to $N_b$, i.e. the index for which
$$ \sum_{j_a+1}^{\tilde{j}_b} \, N_j \;=\; \sum_{i_a+1}^{\tilde{i}_b} \,
N^\prime_{i}
$$

It is clear that $\tau$, correspondingly $\tau^\prime$, is an
element  of $\mathcal{N}(\pi(N))$, respectively of
$\mathcal{N}(\pi^\prime(N))$ of order  $2$  to (a reflection).
Moreover
 $\tau $, correspondingly $\tau^\prime $, does not permute  the summands with indices $\leq 2j_a = \tilde{j}_a=\tilde{i}_a =2 i_a$, and greater then $ \tilde{j}_b$, or $ \tilde{i}_b$ respectively.

 On the other hand, by Proposition
\ref{special partitions} we have $s_N$  not-equivalent classes of
such partitions.

Now, let us fix a frame, i.e. an orthogonal basis $e_1, e_2,\,
\dots\,,e_N$ of $V=\br^N$, Let next  $0 \varsubsetneq V^1
\varsubsetneq V^2 \varsubsetneq \cdots\ \varsubsetneq V^r=V$,
$0\varsubsetneq {V^\prime}^1 \varsubsetneq {V^\prime}^2
\varsubsetneq \cdots\ \varsubsetneq {V^\prime}^s=V $, $\dim
V_j=N_j$, $\dim V^\prime_i =N_i^\prime$ be the partial flags
corresponding to $\pi(N)$ and $\pi^\prime(N)$. Finally let
$$ \ H= O(N_1)\times O(N_2)\times \, \cdots\, \times
O(N_r)\,,\;\;\;\;   H^\prime = O(N^\prime_1)\times O(N^\prime_2)
\times \, \cdots \, \times O(N^\prime_s)
$$
be the maximal orthogonal Borel subgroups corresponding to these
partial flags, i.e. corresponding  to these two partitions. Taking
elements $\tau \in \mathcal{N}(H)\subset  O(N)$, correspondingly
$\tau^\prime \in \mathcal{N}(H^\prime)\subset $, of their
normalizers as described above we get elements which acts
nontrivially only  on the partial flags $V^{j_a+1}\varsubsetneq
V^{j_a+2} \varsubsetneq \, \cdots \, \varsubsetneq V^{\tilde{j}_b}$
and ${V^\prime}^{i_a+1}\varsubsetneq {V^\prime}^{i_a+2}
\varsubsetneq \, \cdots \, \varsubsetneq {V^\prime}^{\tilde{i}_b}$.

We put
$$
\begin{matrix}
V^\prime := \{0\}\times\, \cdots\ \, \times \{0\} \times
\br^{N_{j_a+1}}\times \br^{N_{j_a+2}}\times\, \cdots \, \times
\br^{N_{\tilde{j}_b}}\times \{0\} \times \, \cdots \times \{0\} \cr
{\;\; }\cr  \equiv \{0\}\times\, \cdots\ \, \times \{0\} \times
\br^{N^\prime_{i_a+1}}\times \br^{N^\prime_{i_a+2}}\times  \, \cdots
\, \times  \br^{N^\prime_{\tilde{i}_b}}\times \{0\} \times \, \cdots
 \times \{0\} \simeq \br^{N_{b}-N_a}
\end{matrix}
$$

 Moreover we can take $\tau$  in $\mathcal{N}(H)$, or correspondingly $\tau^\prime$ in  $\mathcal{N}(H^\prime)$, being  of order two as  described above. It is clear that $\tau\in \Upsilon(H)$,
 or  $\tau^\prime$, preserves $V^\prime$ and $\tau \in \Upsilon(H)$,
 or $\tau^\prime\in \Upsilon^\prime(H^\prime)$.

Now, let us consider  the groups $\tilde{H}=\langle H, \Upsilon(H)
\rangle \subset O(N)$ and $\tilde{H^\prime}= \langle H^\prime,
\Upsilon^\prime(H^\prime) \rangle \subset O(N)$.

Finally, take the  corresponding fixed points subspaces $
E^{\tilde{H}}, \;\; E^{\tilde{H\prime}}$, and the subspace
$$\begin{matrix}
 V^\prime := V_{j_a+1} \oplus V_{j_a+2} \oplus \cdots\oplus
V^\prime_{i_b} =V^\prime_{i_a+1} \oplus V^\prime_{i_a+2} \oplus
\cdots\oplus V_{i_b} \simeq  \cr {\;\; }\cr   \br^{N_{j_a+}} \times
\br^{N_{j_a+2}} \times \cdots\times \br^{N_{j_b}}  \simeq
\br^{N^\prime_{i_a+1}} \times \br^{N^\prime_{i_a+2}} \times
\cdots\times \br^{N^\prime_{i_b}} \simeq \br^{N_b-N_a}\,.
\end{matrix} $$

By its construction $V^\prime$ is preserved by $H$ and $H^\prime$
and $\tau$, or $\tau^\prime $ if it is the case. Moreover $\tau$, or
$\tau^\prime$ acts trivially on $V^\prime_\perp$.

 From the above  and Proposition \ref{adaptation of
KrMa} it follows that the group $\langle H, H^\prime\rangle$ acts
transitively on $S(V^\prime)$ and we can apply  Proposition
\ref{linear independence of fixed points}. Consequently
$E^{\tilde{H}} \cap E^{\tilde{H}^\prime} = \{0\}$ if $H$, and
$H^\prime$ are defined by two different partitions.

Let  $\Phi$ be    the functional corresponding to the studied
variational problem with symmetry. Now we can restrict $\Phi$ to
every  subspace $E^{\tilde{H}}$ with $H$ the maximal orthogonal
Borel subgroup as above. By our assumption \ref{Assumption on
functional} $\Phi_{|E^{\tilde{H}}}$ are weak solutions of the
original variational problem by the Palais symmetry principle.

If $u\in E^{\tilde{H}}$ and $u^\prime \in E^{\tilde{H}^\prime}$ then
$u\neq  u^\prime$ as they are in the linearly independent subspaces.
We are left to show that they are geometrically independent, namely
that there is not $g\in  O(N)$ such that $u^\prime(x)= u(gx)$ for
all in $\Omega$. Indeed, if $u(x) \in E^{\tilde{H}} \subset E^H$ and
$u^\prime = u(gx)$, $g\neq e$,  then $u^\prime(x)= u(gx)\in
E^{gHg^{-1}}$. But $gHg^{-1}$ is another maximal orthogonal Borel
subgroup in the same equivalence class as $H$. This implies that the
partition corresponding to $gHg^{-1}$  is the same as that of $H$
which is impossible because we took only one representative $H$  of
each conjugacy class of the maximal Borel subgroups.

Finally, if $u_k^i$, $u_{k^\prime}^i$, $1\leq i \leq s_N$, $k^\prime
> k$,  are two different solutions in one infinite  series $\mathcal{S}_i$
then $\Vert u_{k^\prime}^i\Vert  > \Vert u_k^i \Vert $ by our
assumption \ref{Assumption on functional}. Since  we supposed that
the action of $O(N)$ on $E$ induced by the action of $O(N)$ on
$\Omega$ preserves the norm, it is impossible to have
$u_{k^\prime}^i (x) = u_k^i(gx)$. This shows that these two
solutions are geometrically distinct. \hbx

\begin{rema}\label{second final remark}
Note that Assumption \ref{Assumption on functional} iii) implies
that in each series $\mathcal{S}_i$, $1\leq i \leq s_N$  of Theorem
\ref{final theorem} consists of infinitely many geometrically
distinct functions $\{u^i_k\}$. If we drop out this assumption then
any two different series $\mathcal{S}_i$  of solutions  still
consist geometrically distinct functions, by the same argument as in
the proof of Theorem \ref{final theorem}. Moreover, since they are
given by the variational principle, the cardinality of each series
$\mathcal{S}_i$ is not smaller than the number of critical values of
$\Phi_{|E^{\tilde{G}}}$. Indeed two solutions being in two distinct
critical levels must not be in the same $O(N)$-orbit as $\Phi$ is
$O(N)$-invariant.

\end{rema}
 It is clear that series of solutions found in this way
are not radial and are sign-changing. Indeed $u(\tau( x, y)) =
u(\tau x, y) = -1 u(x,y)$ implies that the nodal  set of $u \in
E^{\tilde{H}}$ contains the fixed point set of reflection $\tau$,
i.e.  a hyperplane $\{(x,y): \tau x = x\}\,.$

In analysis, the questions of finding sign-changing, or
correspondingly  radial solutions for problems posed in $D^N$, or
$\br^N$ is of importance. Obviously a radial solution  is
$SO(N)$-invariant and conversely, so the second question can be
posed as the existence of $SO(N)$-invariant solutions. Of course,
this method gives only some families of the radial and sign-changing
solutions but indicates that this question can be successfully
studied by variational methods as only we are able to defined a
sub-representation of the functional space $E$, which is linearly
independent of $E^{SO(N)}$ is of  the form $E^{G,\rho}$ for a
subgroup $G\subset SO(N)$ and its representation structure $\rho$.
In general, to get not $SO(N)$-invariant, or sign-changing solutions
one can impose some analytical conditions and use subtle analytical
arguments. There  several important works, also very recent,
studying these problems, so we refer only most close to the
analytical problems we have already described \cite{Kaj},
\cite{LoUb}, \cite{MuPaWe}. In the first of quoted papers, the
author showed that for an elliptic problem in $D^N$ and any group
$G\subset SO(N)$ such that $G$ does not act transitively on
$S(\br^N)$ there exists infinitely many solutions which are
$G$-invariant but not $SO(N)$-invariant (not radial). He did it by
proving that the distribution of critical values corresponding  to
the  functional restricted to $E^{SO(N)}$ is smaller than the
distribution of critical values which correspond to the restriction
of functional to $E^G$. In the second paper the authors studied a
similar  problem as that of \cite{BaWi} for $N=5$ which dimension is
not covered by the approach used in \cite{BaWi} and studied here for
the obvious combinatorial reason. But they impose a condition of the
form nonlinear term which let them to reduce the studied  problem to
a problem  in $\br^4$. Finally, in \cite{MuPaWe} for    a problem
like in \cite{BaWi} and every  $N\geq2$ but with some additional
assumption on the nonlinear part it is shown that for any $k\geq 7$
there exist infinitely many $D_k \times O(N-2)$-invariant solutions,
where $D_k\subset O(2)$ is the dihedral group of $2k$ elements.
Moreover these solutions are not $O(2)\times O(N-2)$-invariant, e.g.
for $N=4$ give a nice example of not $O(N)$-invariant (radial)
solutions which are different than those obtained by the method of
\cite{BaWi}, i.e. by the scheme studied here.

\subsection{Spaces of $\rho$-interwinding functions}\label{rho-interwinding subspaces}
It is possible to define these functional subspaces associated with
a partition $\pi(N)$  in terms of representations of the Weyl group
$\mathcal{W}(H)$, where $H$ is the unique maximal orthogonal Borel
subgroup defined by $\pi(N)$. Then  it would be seen that the nodal
set contains union of fixed points of all elements of order two in
$\mathcal{W}(H)= \mathcal{W}(\pi(N))$, i.e. the union of  fixed
points of all transposition, or equivalently reflections, in
$\mathcal{W}(\pi(N)) = \mathcal{W}(H)$. But we have consider these
reflections, and their fixed points, only in subspaces on which they
are not trivial.

Let  $\{V_j\}_1^r=\{\br^{N_1}, \br^{N_2}, \dots, \br^{N_r}\}$ be an
orthogonal partial flag in $\br^N$ with $\dim V_j =N_j$, or
equivalently the corresponding  maximal orthogonal  Borel subgroup
uniquely determined the partition $\pi$ of $N= \sum_{j=1}^r\, N_j$.

By Proposition \ref{Weyl group of flag} and Corollary  \ref{Weyl
group of subgroup}  the normalizer $\mathcal{N}(\{V_j\}_1^r)$ in
$O(N)$, or equivalently the normalizer $H$ in  $O(N)$   consists of
the maximal orthogonal Borel subgroup $H=O(N_1)\times O(N_2)\times\,
\cdots\, \times O(N_r)$  and all permutations
$$\sigma_\pi \in \mathfrak{S}_\pi =  \mathfrak{S}(\phi_\pi(n_1))
\times \mathfrak{S}(\phi_\pi(n_2)) \times
    \cdots\, \times
\mathfrak{S}(\phi_\pi(n_q)\,.$$

Moreover  the Weyl group of $H$ is a product of $q$ permutation
groups
$$ \mathcal{W}(H)= \mathfrak{S}_\pi =  \mathfrak{S}(\phi_\pi(n_1)) \times \mathfrak{S}(\phi_\pi(n_2)
\times
    \cdots\, \times
\mathfrak{S}(\phi_\pi(n_q))$$

Let $\rho^1:\mathfrak{S}(n) \to \{-1,1\} = O(1)$ be the unique
non-trivial one-dimensional representation  of the permutation group
$\mathfrak{S}(n)$ given by $\sigma \mapsto {\rm sign}  \sigma$.
Correspondingly, let $\rho^0:\mathfrak{S}(n) \to \{-1,1\} = O(1)$ be
the trivial one-dimensional representation  of the permutation group
$\mathfrak{S}(n)$ given by $\sigma \mapsto 1$ for all $\sigma \in
\mathfrak{S}(n)$.

Moreover, since $\mathcal{W}(H)= \mathfrak{S}(\phi_\pi(n_1)) \times
\mathfrak{S}(\phi_\pi(n_2) \times
    \cdots\, \times
\mathfrak{S}(\phi_\pi(n_q))$ is the product of group, every
one-dimensional orthogonal representation  of $\mathcal{W}(H)$, i.e
.every homomorphism $\rho$ from $\mathcal{W}(H)$ to $\bz_2= O(1)$ is
of the form
$$ \rho(\sigma_1, \, \sigma_2, \,\dots\, \sigma_p)=
\rho_1^{\delta_1}(\sigma_1) \cdot \rho_2^{\delta_2}(\sigma_2) \,
\cdots\, \rho_p^{\delta_p}(\sigma_p)$$ where $\delta_i \in \{ 0,1\}
$, and $\rho_i^{\delta_i}:\mathfrak{S}(\phi_\pi(n_i)) \to \{-1,1\}$
is either the the unique nontrivial homomorphism, or correspondingly
the trivial homomorphism from $\mathfrak{S}(\phi_\pi(n_i))$ to
$O(1)$ depending whether $\delta_i$ is equal to $1$, or to $0$
respectively.  In other words every such representation is
determined by the sequence $(\delta_1,\, \delta_2, \, \cdots,\,,
\delta_p)$, e.g. $\rho$ is trivial if and only if this sequence
consists of zeros only.

Furthermore, every representation $\rho$  of $\mathcal{W}(H)$
composed with the natural projection $\mathcal{N}(H) \to
\mathcal{N}(H)/H= \mathcal{W}(H)$ defines a representation of
$\mathcal{N}(H)$. Observe that  $\rho_{|H} \equiv 1$ for every such
$\rho$.

\begin{defi}\label{special rho-interwinding functions}
Let $ H = O(N_1)\times  O(N_2) \times \, \cdots\, \times O(N_r)$ be
an maximal orthogonal Borel subgroup with  nontrivial the Weyl group
$\mathcal{W}(H)= \mathfrak{S}(\phi_\pi(n_1)) \times
\mathfrak{S}(\phi_\pi(n_2) \times
    \cdots\, \times
\mathfrak{S}(\phi_\pi(n_q))$. Let $\rho= \rho_1^{\delta_1} \,
\rho_2^{\delta_2}\, \cdots\, \rho_q^{\delta_q}$ be a nontrivial
one-dimensional representation of $\mathcal{N}(H)$. Then in  every
functional space $E$ with domain $\Omega = \br^N$, or $D^N$ we can
define a linear action of $\mathcal{N}(H)$ (a representation
structure) by the formula
$$ \begin{matrix}  g u(x):= \rho(g) u(gx),\;\;  {\text{or equivalently for}} \cr  \; g= \sigma \,h,
\; \sigma \in \mathcal{W}(H)\;\; {\text{ and}}\; \; h\in H\;\;
{\text{by}} \;\; (\sigma, h ) u(x):= \rho(\sigma) u(\sigma h x)\,,
\end{matrix}
$$
 where $\sigma =
(\sigma_1, \sigma_2,\, \dots \,,\sigma_p)$, and the formula does not
depend on a representation of $g$ as a pair $(\sigma, h)$ since
$\rho_{|H} \equiv 1$.

Finally to a pair $(H,\rho)$ as above we assign the fixed point
space of this action denoted by $E^{(H,\rho)}$ and called  (as in
\cite{BrClMa}) the space of {\underline{$\rho$-interwinding
functions}}.
\end{defi}
\begin{rema}
Note that the spaces $E^{\tilde{H}}$ discussed previously  are
special cases of spaces defined in Definition \ref{special
rho-interwinding functions}, i.e. every space $E^{\tilde{H}}$ is of
this form.

Our task  was to find a  family  $(H_i,\rho^i)$, $1\leq i \leq s_N$
as above such that $E^{(H_,\rho^i)} \cap E^{(H_j,\rho^j)} =\{0\}$ if
$i\neq j$ and with possibly large $s_N$.
\end{rema}

Observe that each $\mathfrak{S}(\phi_\pi(n_1))$ is generated by
transpositions $\tau$, geometrically reflections, so
$\mathcal{W}(H)= \mathfrak{S}(\phi_\pi(n_1)) \times
\mathfrak{S}(\phi_\pi(n_2) \times
    \cdots\, \times
\mathfrak{S}(\phi_\pi(n_q))$ is generated by the compositions of
transpositions. As $ u(\tau x) = -1 u(x) $ provides $\tau\in W(H)$
and $\rho(\tau)=-1$. The latter implies that $u(x)= 0  $ if $\tau x
=x $ is a fixed point of $\tau$. As a consequence we get the
following.
\begin{coro}\label{nodal set}
The zero set (the nodal set if $u$ is a solution) of every $u\in
E^{(H,\rho)}$ with $\rho=
\rho_1^{\delta_1}\rho_2^{\delta_2}\,\cdots\, \, \rho_q^{\delta_q}$
contains
$$ {\underset{\delta_i =1}\bigcup} \, \mathcal{H}_i $$
where $\mathcal{H}_i$ is a union of hyperplanes being fixed points
of transposition $\tau \in \mathfrak{S}(\phi_\pi(n_i))$ interpreted
as reflections.

In particular, if $H= O(2)\times O(2)\times \, \cdots\, O(2) \subset
O(2N)$, then $\mathcal{W}(H)= \mathfrak{S}(N)$ is equal to
$W(\mathbb{T}^N) \subset SO(2N)$. Consequently for every $u\in
E^{(H,\rho)}$ its zero set contains the union of walls of the Weyl
chambers of canonical representation of $O(N)$ in $\br^N$.
\end{coro}
\begin{rema}
It is worth of pointing out that the functional  subspaces  of
$\rho$ -interwinding functions as defined in Definition \ref{special
rho-interwinding functions} can be define in a context of  action of
any Coxeter group, but we do not know any application, especially
that the part which shows that some of them are orthogonal has not a
direct analog.
\end{rema}

\section{Supplementary information}\label{Appendix}

\subsection{Borel groups} In 1950 A. Borel \cite{Bor1} gave a complete
classification of all groups acting transitively and effectively on
the sphere completing earlier results of D. Montgomery and H.
Samelson \cite{MoSa}.

\begin{teo}[\rm A. Borel]\label{Borel theorem}
If a  connected compact group $G$ of linear transformations acts
effectively transitively on the sphere $S(\mathbb{R}^N)$ then
$S(\mathbb{R}^N)$ is $G$-homeomorphic to the homogenous space $G/H$
where the pair $(G,H)$ is one of the listed below.

\begin{itemize}
\item[i)] {If $N$ is odd, or equivalently $N-1$ is even, then $(G,H)\simeq (SO(N), SO(N-1))$, or $(G_2,
SU(3))$ in the case when $N=7$;}
\item[ii)] {If $N=2s$, or equivalently $N-1=2s-1$ and $s$ is odd, then   $(G,H)\simeq (SO(N), SO(N-1))$ or $(SU(s),
SU(s-1))$;}

\item[iii)] {If $N=2s$, or equivalently  $N-1 = 2s-1$ and $s$ is even,  then $(G,H) \simeq (SO(N), SO(N-1))$,
$(SU(s), SU(s-1))$, $(Sp(s/2), (Sp(s/2-1))$,}

\item[iv)]  {$(Spin(9), Spin (7))$ in case $n=15 $, or $(Spin (7), G_2)$ in
case $n=7$.}
\end{itemize}
\end{teo}

Here $G_2 \subset GL(7, \br)$ denotes the automorphism group of the
octonion algebra, i.e. the subgroup of $GL(7,\br)$ of  that
preserves the non-degenerate 3-form

$$ dx^{123} + dx^{235} + dx^{346} +dx^{450} + dx^{561} +dx^{602}
+dx^{013} $$ (invariant under the cyclic permutation (0123456)) with
$dx^{ijk}$ denoting $dx^j \wedge dx^j \wedge dx^k$ in variables
$x^i, \, 0\leq i \leq 6$ of $\br^7$.

\begin{rema} It is worth to pointing out that original formulation
of the Borel theorem is stronger, i.e. it has a weaker  supposition
that $G$ is a  compact connected  Lie group of transformations of a
homotopy sphere acting effectively and transitively.
\end{rema}

The Borel theorem says, roughly speaking,  that only $SO(N)$ or in
few cases its classical linear subgroups are  only connected groups
that act transitively on $S(\br^N)$. The latter happens if $\br^N$
has an extra structure: complex, quaternionic, spinor, or octonion.

We end this subsection with a statement which is generalization of
Proposition \ref{adaptation of KrMa} and a positive answer to a
question posed in \cite[Remark 4.1]{KrMa} in a stronger form.

Once more, let $\pi_H = \{ N^H_1 \leq N^H_2\leq N^H_3\leq \,
\cdots\, \leq N^H_r\}$, $\pi_K= \{N_1^K \leq N^K_2\leq \,
\cdots\,\leq N^K_s\}$, \newline  $\sum_1^r N_j^H = \sum_1^s N^K_i
=N$ be two different not decreasing partitions  of $N$ corresponding
to the maximal orthogonal Borel  subgroups $H=O(N^H_1)\times
O(N^H_2) \times\, \cdots\, \times O(N^H_r)$, \newline
$K=O(N^K_1\times O(N^K_2) \times\, \cdots\, \times (N^K_s)$, or
maximal orthogonal connected Borel subgroups $H=SO(N^H_1)\times
SO(N^H_2) \times\, \cdots\, \times SO(N^H_r)$, $K=SO(N^K_1)\times
SO(N^K_2)\, \times\, \cdots\, \times SO(N^K_s)$ respectively.
Remind, that we say that $\pi_H$ and $\pi_K$ contain common
partition of $N^\prime <N$ if there exist $1\leq a\leq r$, $1\leq
b\leq s$ such that $\sum_{j=1}^a N^H_j = \sum_{i=1}^b
N^{K}_j=N^\prime$. We begin with a proposition which is a positive
answer to the question posed in  \cite{KrMa}.

\begin{prop}\label{KrMrTheorem}
Let $K, \, H \subset O(N)$ be two maximal orthogonal Borel subgroups
as above, respectively   $K, \, H \subset O(N)$ be two maximal
orthogonal connected Borel subgroups,  corresponding to two
different partitions $\pi_H$ and $\pi_K$, of $N$ such that
$N^H_i\geq 2 $ and $ N^K_j\geq 2$ for all  $1\leq i\leq r$, $1 \leq
j \leq s$.

Then $\pi_H$ and $\pi_K$ do not contain a common partition of
$N^\prime <N$ if and only if $ \langle H, \,K\rangle = O(N)$,  or
respectively $ \langle H, \,K\rangle = SO(N)$ when $H, \,K$ are
connected.
\end{prop}

\noindent{\sc{Proof.}} We will show the connected case only.

 The
part "if" is a direct consequence of \cite[Lemma 4.1]{KrMa} shown in
\cite[Remark 4.1]{KrMa}. To prove the part "only if" observe that $
\langle H, \,K\rangle$ acts transitively on $S(\mathbb{R}^H)$ as
follows from  \cite[Lemma 4.1]{KrMa} (cf. Proposition
\ref{adaptation of KrMa}). Consequently it is one of the groups
listed in the statement of Borel theorem \ref{Borel theorem}.

 By the
dimension assumption $N^H_1\geq 2$, $N^K_1\geq 2$. Next, since
$\pi_H$ and $\pi_K$ are different and do not contain a common
partition of $N^\prime < N$, either $N^H_1\geq 3$ or $N^K_1\geq 3$.
Consequently the group $SO(3)$ is contained  at least in one of the
subgroups $H$ or $K$. Suppose that $N^H_1\geq 3$.  Is enough to find
an orthogonal map $ A \in SO(3) \subset SO(N^H_1)$ and its extension
$\tilde{A}$ to $H$ such that $A \notin G$ for all $G \subset SO(N)$
listed in Theorem \ref{Borel theorem}. Consider $A$ given by the
matrix
$$ \begin{bmatrix} -1 & 0 & 0 \cr 0 & 1 & 0 \cr 0 & 0 & - 1
\end{bmatrix} $$ and extended to a linear map of $\br^{N_1}$ by the
identity on remaining coordinates.

Since $A\in SO(\br^{N_1})$, the element $ \tilde{A}= A\times \{e\}
\,\times \{e\}\, \times \cdots\, \times \,\{e\}$ belongs to $H$,
thus belongs to $\langle H, \, K\rangle $. On the other hand
$\tilde{A} $ does not belong to any of proper subgroup $G$ of
$SO(N)$ listed in the Borel theorem \ref{Borel theorem}. This shows
that $\langle H, \,K\rangle =SO(N)$ which completes the proof. \hbx

Finally, we are able to describe a subgroup generated by two maximal
orthogonal subgroups $H, \,K \subset O(N) $ associated with two
partitions $\pi_H$, $\pi_K$ of $N$. Reasoning as in before
Proposition \ref{adaptation of KrMa} we split the sequences  of (non
decreasing) partitions $\{N^H_1, \, N^H_2, \, \cdots \,,  N^H_r\}$,
$\{K^H_1, \, N^K_2, \, \cdots \,, N^K_s\}$ subsequences following
intervals of indices for which these partitions are equal or not
comparable. More precisely, let $1\leq j_{a_1} < j_{b_1} \leq
j_{a_2} < j_{b_2} \, \dots\, \leq j_{a_k}$ and $1\leq i_{a_1} <
i_{b_1} \leq i_{a_2} < i_{b_2} \, \dots\, \leq i_{a_k}$ be indices
such that $j_{a_1}= i_{a_1}$ and $N^H_j = N^K_i$ for all $j, \,
i\leq j_{a_1}$,  in the next  interval of indices  $j_{a_1} <
j_{{a_1}+1} <\dots j_{b_1}$, $i_{a_1} < i_{{a_1}+1} <\dots i_{b_1}$
be such that $\pi_H$ and $\pi_K$ are not commensurable, i.e. do not
exist $j_{a_1} <j_{c_1} < j_{b_1}$, $i_{a_1} < i_c < i_{b_1}$ such
that $\sum_{j=j_{a_1}}^{j_{c_1}} N^H_j =  \sum_{i=i_{a_1}}^{i_{c_1}}
N^K_i$ but $$\sum_{j=j_{a_1}}^{j_{b_1}} N^H_j =
\sum_{i=i_{a_1}}^{i_{b_1}} N^K_i = N_{b_1} -N_{a_1}\,,$$ where
$N_{a_1} = \sum_{j=1}^{j_{a_1}} N^H_j$, $N_{b_1} =
\sum_{j=1}^{j_{b_1}} N^H_j = \sum_{i=1}^{i_{b_1}} N^H_{i_{b_1}} $.
In the next intervals $j_{b_1} \leq j_{a_2}$, and $i_{b_1} \leq
i_{a_2}$ the partitions are equal and so on alternately. Put
$N_{b_j, a_{j}}:= N_{b_i} - N_{a_i} = N_{b_j}- N_{a_j}$.

\begin{teo}\label{group generated by aprtitions} Let  $\pi_H=\{N^H_1, \, N^H_2, \, \cdots \, ,N^H_r\}$,
and $\pi_K=\{K^H_1, \, N^K_2, \, \cdots \,, N^K_s\}$ be two
partitions of $N$. Let next  $H=O(N^H_1)\times O(N^H_2) \times\,
\cdots\, \times O(N^H_r)$, $K=O(N^K_1)\times O(N^K_2) \times\,
\cdots\, \times  O(N^K_s)$ be associated with them maximal
orthogonal Borel subgroups of $O(N)$, correspondingly
$H=SO(N^H_1)\times SO(N^H_2) \times\, \cdots\, \times SO(N^H_r)$,
$K=SO(N^K_1)\times SO(N^K_2) \times\, \cdots\, \times SO(N^K_s)$
maximal connected orthogonal Borel subgroups.

With the above notation, the group  $\langle H, \, K\rangle $ is
equal to
$$ O(N^H_1)\times O(N^H_2)\times \dots \times O(N_{j_{a_1}})\times
O(N_{b_1, a_1}) \times O(N^H_{j_{b_1}+1}) \times O(N_{j_{a_2}})
\times \dots O(N^H_{j_{a_2}}) \times \, \dots\, \times
O(N^H_{b_k,a_m})$$ with the convention that if
$j_{b_q}=j_{a_{q+1}}$, or correspondingly $b_k=a_m$ then the
corresponding factor $O(N^H_{j_{b_q}+1}) \times O(N_{j_{a_{q+1}}})
$, or $O(N^H_{b_k,a_m})$ respectively,  is equal to $\{e\}$.
Moreover the thesis holds for the connected case with a
corresponding formulation.
\end{teo}
\noindent{\sc{Proof.}} The statement is a direct consequence of
Proposition \ref{KrMrTheorem} applied  consecutively to the
partitions $N^H_{j_{a_q}}, N^H_{j_{a_q} +1}, \dots, N^H_{j_{b_q}}$
and $N^H_{i_{a_q}}, N^H_{i_{a_q} +1}, \dots, N^H_{i_{b_q}}$. In the
intervals in which the partitions are equal the factors of $\langle
H,\, K\rangle $  are the same as in the original groups $H$ and $K$.
\hbx

\subsection{Partitions and their properties}

The functions  $P(N)$ and $Q(N)$ have been studied by several
mathematicians for hundreds years, so let us give only  references
to the survey articles from MathWorld--A Wolfram Web Resource
(\cite{Wei1, Wei2}) where one can find as well basic facts as an
expanded bibliography.

\begin{rema}\label{computations of R(N)}
It is worth of pointing out that the function $R(N)$ is strictly
monotonic for $N\geq 4$, and the first its values, for $N=1, \,
\dots,\, 10 $ are equal to
$$  
 { {0,
1, 1, 3, 4, 7, 10, 16, 22, 32}}
$$
\end{rema}

\begin{rema}
Note that the functions $ P(N; 1), \, Q(N; 1)\,, R(N; 1)$ are not
expressed  directly by the classical number theory functions $
P(N,k), \, Q(N, k)$ of the number of partitions of $N$ into summands
each of them is smaller or equal to $k$.
\end{rema}

\begin{lema}\label{formula for P(N; 1)}
For the functions $ \,P(N; 1)\,,Q(N; 1)$ of Definition
\ref{partition without 1} we have
$$ {\text{For}}\;\;N\geq 2\;\;  \,P(N; 1)\, = P(N)- P(N-1)$$
$$ {\text{For}}\;\;N\geq 2\;\;  \,Q(N; 1)\,= Q(N) - Q(N-1; 1) $$ therefore $$\; Q(N; 1) =
Q(N) - Q(N-1) + Q(N-2)+\, \cdots \, +(-1)^{N-2} Q(2)\,.$$
Consequently
$$\,R(N; 1)\,= P(N)-  P(N-1) \;-\; \Big(\sum_{i=0}^{N-2} (-1)^i
Q(N-i)\Big)\,.$$
\end{lema}
\noindent{\sc Proof.} Since $P(N)$ measures the number of partitions
of $N$ where the order of addends is not considered significant, we
can assume that in every such partition $N_1, N_2, \, \dots\,,N_r$
we have $N_1\leq N_2\leq \, \dots\,\leq N_{r-1}\leq N_r$. If $N_1=1$
then $N_2,\, N_3,\,\dots\,,N_r$ gives a partition of $N-1$ of length
$r-1$. Conversely, if $N_2,\, N_3,\,\dots\,,N_r$ is a partition of
$N-1$ of length $r-1$ then $1\leq N_2\leq \, \dots\,\leq N_{r-1}\leq
N_r$ is a partition of $N$, because $N_2 \geq 1$. Consequently the
only partitions of $N$ for which $N_1\geq 2$ are exactly these which
are not constructed by the above procedure, which shows that $P(N;
1)\,= P(N) - P(N-1)$.

Next, let $N_1< N_2< \, \dots\,< N_{r-1}< N_r$ is a partition of $N$
consisting of distinct summands. If $N_1=1$ then $N_2,\, N_3,\,
\dots\,, N_r$ is a partition of $N$ of length $N-1$ consisting of
distinct summands and not containing $1$. Conversely, if  $2\leq
N_2< \, \dots\,< N_{r-1}< N_r$ is  a partition of $N-1$ of length
$r-1$ consisting of distinct summands and not containing $1$, then
adding $1$ as first summand we get a partition  $ 1 < N_2< N_3<\,
\dots\, <N_r$ of $N$ of length $r$. This shows that $Q(N; 1)\,= Q(N)
- Q(N-1; 1)$. The last equality follows from the latter by applying
it $N-1$ times. \hbx

\vskip 0.3cm
 For  a convenience of  the  reader we present the
values of $P(N),\, Q(N),\, P(N; 1), \, Q(N; 1)$, and $R(N; 1)$ for
first $10$ values of $N$
\begin{abc}\label{table1}\rm\hskip 2cm
\begin{tabular}{|c|c|c|c|c|c|c|}
  \hline
  N & $P(N)$ & $Q(N)$ & $R(N)$ & $P(N; 1)$ & $Q(N; 1)$ & $R(N; 1)$ \\
  \hline
  1 &  1 & 1 & 0 & 0 & 0 & 0 \\ \hline
  2 & 2 & 1 & 1 & 1 & 1 & 0 \\ \hline
  3 & 3 & 2 & 1 & 1 & 1 & 0 \\ \hline
  4 & 5 & 2 & 3 & 2 & 1 & 1 \\ \hline
  5 & 7 & 3 &  4 & 2 & 2 & 0 \\ \hline
  6 & 11 & 4 & 7 & 4 & 2 & 2 \\ \hline
  7 & 15 & 5 & 10 & 4 & 2 & 2 \\ \hline
  8 & 22 & 6 & 16 & 7 & 4 & 3 \\ \hline
  9 &  30 & 8 & 22 & 8 & 5 & 3 \\ \hline
  10 &  42 & 10 & 32 & 12 & 5 & 7 \\ \hline
  \hline
\end{tabular}

\end{abc}

\vskip 0.5cm The values of $P(N)$ for first $49$ natural
numbers:{\footnotesize{\small{
\newline

\noindent $N=1, \,P(N)=1; \; N=2,\, P(N)=2; \; N=3,\, P(N)=2;\;
N=4,\,P(N)= 5;\; N=5,\, P(N)=7; $  $\; N=6,\, P(N)=11;\; N=7,\,
P(N)=15;\;N=8,\, P(N)=22;\;N=9, P(N)=30;\;N=10, P(N)=42;$  $\;N=11,
P(N)=56;\;N= 12,P(N)=77;\;N=13,\, P(N)=101;\;N=14,\,
P(N)=135;\;N=15,\, P(N)=176;$  $\; N=16,\, P(N)=231;\; N=17,
\,P(N)=297;\; N=18, \,P(N)=385; \; N=19,\, P(N)= 490;\; N=20,\,
P(N)= 627;$ $\; N=21, P(N)=792;\; N=22,\, P(N)= 1002;\; N=23,\,
P(N)= 1255; \; N=24,\,P(N)= 1575; \; N= 25,\, P(N)=1958;$ $\;
N=26,\,P(N)= 2436; \; N=27,\, P(N)=3010; \; N=28,\, P(N)=3718; \;
N=29, \,P(N)=4565; \; N=30, \,P(N)=5604; $  $\; N=31,\, P(N)=6842;
\; N= 32, \, P(N)=8349; \; N=33, \, P(N)= 10143; \; N=34, \, P(N)=
12310; \; N=35,\,P(N)= 14883;$
 $\;N=36,
\, P(N)= 17977;\;N=37,\, P(N)= 21637; \;N= 38,\, P(N)= 26015;\; \;
N= 39,\, P(N)= 31185;\;  N=40,\, P(N)= 37338;$
 $\; N= 41,\, P(N)=
44583;\; N= 42,\, P(N)= 53174; \; N= 43, \,P(N)= 63261; \; N= 44,\,
P(N)= 75175; \; N= 45,\, P(N)= 89134;$  $\; N=46,\, P(N)=105558; \;
N=47,\, P(N)= 124754; \; N= 48,\, P(N)= 147273; \; N=49, \,P(N)=
173525.$}}}

    \end{document}